\documentclass{amsart}[11pt]
\usepackage{amsmath, amssymb, amsthm, graphics}
\theoremstyle{plain}
\newtheorem{prop}{Proposition}[section]
\newtheorem{coro}[prop]{Corollary}
\newtheorem{lemm}[prop]{Lemma}
\newtheorem{ques}[prop]{Question}
\newtheorem{thrm}[prop]{Theorem}

\theoremstyle{definition}
\newtheorem*{defi}{Definition}
\newtheorem*{nota}{Notation}
\newtheorem{exam}[prop]{Example}
\newtheorem{rema}[prop]{Remark}

\numberwithin{equation}{section}

\def\Reff#1; #2; #3; #4; #5; #6; #7\par{%
\bibitem{#1} #2, {\it #3}, #4 {\bf #5} (#6) #7}

\def\Ref#1; #2; #3; #4\par{%
\bibitem{#1} #2, {\it #3}, #4}

\let\a=\alpha

\let\b=\beta
\def\bb(#1, #2){\langle #1, #2\rangle}		
\def\BKL{{B_{_{\!B\!K\!L}}}}
\def\BKL{B\!K\!L}

\def\CC#1{C_{#1}}
\def\CCC#1{\undertilde 3 C_{#1}}
\def\c[#1]{[#1]} 
\def\cc[#1,#2]{[#1, #2]} 
\def\ccc[#1,#2,#3]{[#1, #2, #3]} 
\def\cccc[#1,#2,#3,#4]{[#1, #2, #3, #4]} 
\def\cl#1{\overline{\vrule width 0pt height 5pt#1}}
\def\cl#1{\overline{#1}}
\def\comp{\mathbin{\scriptscriptstyle \circ}}

\let\d=\delta
\let\D=\Delta
\def\DD#1{\hbox{d}_{#1}}
\def\DDD#1{\undertilde3{\hbox{d}}_{#1}}
\def\dd#1{\partial_{#1}}
\def\ddd#1{\undertilde3\partial_{#1}}
\def\DivD{{\mathcal D}_{\!_\Delta}}
\def\div{\mathrel{\scriptstyle\sqsubset}}
\def\dive{\mathrel{\scriptstyle\sqsubseteq}}

\let\e=\varepsilon
\def\ec{\c[\emptyset]}		

\def\first#1{{#1}_{(1)}}

\let\g=\gamma
\def\Gr(#1; #2){\langle #1 \, ; \, #2 \rangle}

\let\Hat=\widehat

\def\id{{\rm id}}
\def\ie{{\it i.e.}}
\def\ii{^{-1}}
\def\Im{{\rm Im}}

\def\Ker{{\rm Ker}}

\def\lc(#1, #2){{#1}_{\! / \!{ #2}}}
\def\lcm{\vee}
\def\LCM#1{^\ulcorner\!\! #1\hspace{-1pt}_\lrcorner}
\def\slashstar{\smash{/^{\!^{^*}}}}
\def\lR(#1, #2){{#1}_{\! \slashstar \!\!\! {#2}}}

\def\slashstar{\smash{/^{{\!*}}}}
\def\lR(#1, #2){{#1}_{\! \slashstar \! {#2}}}
\def\lo(#1, #2){\lR(#1, #2)}
\def\lRD(#1, #2){\lR(#1, #2)}
\def\lRo(#1, #2){\lR(#1, #2)}

\def\md#1{\hbox{mindiv}(#1)}		
\def\Mon(#1; #2){\langle #1 \, ; \, #2
\rangle^{\scriptscriptstyle +}}
\def\mul{\mathrel{\scriptstyle\sqsupset}}
\def\mule{\mathrel{\scriptstyle\sqsupseteq}}

\def\NF#1{\undertilde 5 {\hbox{\rm\small N\!F}}(#1)}
\def\NFD#1{\mbox{\rm\small N\!F}(#1)}

\let\o=\omega
\let\op=\cdot

\def\PH#1{(H_{#1})}
\let\pp=\ldots
\def\PP#1{(P_{#1})}
\def\PQ#1{(Q_{#1})}

\def\rc(#1, #2){_{#1 \! \backslash} \! #2}
\def\RX{{\mathcal R}_{\!{}_\XX}}  
\def\red#1{r_{#1}}
\def\resp{resp\hbox{. }}
\def\rr{\rho} 		
\def\RO{{\mathcal R}_{\!{}_<}}  
\def\RRR{{\mathcal R}}  

\let\s=\sigma

\def\ss#1{s_{#1}}  
\def\sss#1{\undertilde 2 s_{#1}}  

\def\trX{\triangleright_\XX}
\def\TTO#1#2{\mathrel{\raise-3pt\hbox{$\buildrel 
{\overset{\textstyle #1}{\textstyle\longrightarrow}} 
\over {\underset{\textstyle #2}{\longleftarrow}}$}}}

\def\undertilde#1#2{\hskip#1pt
\smash{\raise-8pt\hbox{$\tilde{}$}}
\hskip-#1pt#2}

\def\XX{{\mathcal X}}
\def\XXX#1{\XX^{[\undertilde{2.5}{#1}]}}

\def\ZZ{{\mathbf Z}}

\begin{document}

\author{Patrick DEHORNOY}
\address{Laboratoire SDAD, Math\'ematiques\\
Universit\'e Campus II, BP 5186, 14032 Caen, France}
\email{dehornoy@math.unicaen.fr}
\urladdr{//www.math.unicaen.fr/\!\!\!\hbox{$\sim$}dehornoy}

\author{Yves LAFONT}
\address{IML, 163, avenue de Luminy, 13288
Marseille cedex~9, France}
\email{lafont@iml.univ-mrs.fr}

\title{HOMOLOGY OF GAUSSIAN GROUPS}



\begin{abstract}
We describe new combinatorial methods for
constructing an explicit free resolution of~$\ZZ$ by
$\ZZ G$-modules when $G$ is a group of fractions of a
monoid where enough least common multiples exist
(``locally Gaussian monoid''), and, therefore,
for computing the homology of~$G$. Our
constructions apply in particular to all Artin
groups of finite Coxeter type, so, as a
corollary, they give new ways of computing
the homology of these groups.
\end{abstract}

\maketitle

\section*{Introduction}

The (co)homology of Artin's braid groups~$B_n$ has
been computed by methods of differential geometry
and algebraic topology in the beginning of the 1970's
\cite{Arn, Arp, Fuk, Coh}, and the results have then
been extended to Artin groups of finite Coxeter type
\cite{Bri, Gor, Vai}, see also~\cite{Coi, CSS, Sal, Ver}. A
purely algebraic and combinatorial approach was
developed by C.~Squier in his unpublished PhD thesis
of~1980---see~\cite{Squ}---relying both on the fact that
these groups are groups of fractions of monoids
admitting least common multiples and on the
particular form of the Coxeter relations involved in
their standard presentation.

On the other hand, it has been observed in recent
years that most of the algebraic results established
for the braid groups and, more generally, the
Artin groups of finite Coxeter type (``spherical Artin
groups'') by Garside, Brieskorn, Saito, Adyan, Thurston
among others, extend to a wider class of so-called
Garside groups. A Gaussian group is defined to be the
group of fractions of a monoid in which left and right
division make a well founded lattice, \ie, in which we
have a good theory of least common multiples, and a
Garside group is a Gaussian group that satisfies an
additional finiteness condition analogous to
sphericality (see the precise definition in
Section~\ref{S:Gaus} below). In some sense, such an
extension is natural, as the role of least common
multiples (lcm's for short) in some associated monoid
had already been emphasized and proved to be crucial
in the study of the braid groups, in particular in the
solution of the conjugacy problem by Garside
\cite{Gar} and the construction of an automatic
structure by Thurston \cite{Thu}, see also \cite{Eps,
Cha, Chb}. However, the family of Garside groups
includes new groups defined by relations
quite different from Coxeter relations, such as $\langle
a, b, c, \pp ; a^p = b^q = c^r = \pp \rangle$, $\langle a,
b, c ; abc = bca = cab \rangle$, or $\langle a, b ; ababa
= b^2 \rangle$---see \cite{Pic} for many
examples---and, even if the fundamental
K\"urzungslemma of~\cite{BrS} remains valid in all
Gaussian monoids, many technical results about
spherical Artin groups fail for general Gaussian groups,
typically all results relying on the symmetry of the
Coxeter relations, like the preservation of the length by
the relations or the result that the fundamental
element~$\D$ is squarefree. Thus, the extension from
spherical Artin groups to general Gaussian groups or,
at least, Garside groups is not trivial, and, in most
cases, it requires finding new arguments: see
\cite{Dfx} for the existence of a quadratic isoperimetric
inequality, \cite{Dfz} for torsion freeness,
\cite{Dgk} for the existence of a bi-automatic structure,
\cite{Pid} for the existence of a decomposition into a
crossed product of groups with a monogenic center,
\cite{Sib} for the decidability of the existence of
roots.

According to this program, it is natural to look for a
possible extension of Squier's approach to arbitrary
Gaussian groups (or to even more general groups).
Such an idea is already present in Squier's paper, whose
first part addresses general groups and monoids which
are essentially the Gaussian groups we shall consider
here. However, in the second part of his paper, he can
complete the construction only in the special case of
Artin groups. Roughly  speaking, what we do in the
current paper is to develop new methods so as to
achieve the general program sketched in the first part
of~\cite{Squ}.

As in~\cite{Squ}, we observe that the homology of a
group of fractions coincides with that of the involved
monoid, so our aim will be to construct a resolution
of~$\ZZ$ by free $\ZZ M$-modules when $M$ is a
monoid with good lcm properties. In the spirit of the
standard resolution, we start with the natural idea of
constructing an explicit simplicial complex where the
$n$~cells correspond to $n$-tuples of elements
$(\a_1, \pp, \a_n)$ of~$M$, but, in order to obtain
smaller (finite type) modules, we assume in addition
that the~$\a_i$'s are taken in some fixed set of
generators of~$M$. The idea, which is already present
in~\cite{Squ} even if not stated explicitly, is that the
cell $\ccc[\a_1 , \pp , \a_n]$ represents in some sense
the computation of the left lcm of~$\a_1$, \pp,
$\a_n$. The core of the problem is to define the
boundary of such a cell and to construct a contracting
homotopy. Here Squier uses a trick that allows him to
avoid addressing the question directly. Indeed, he first
defines by purely syntactical means a top degree
approximation of the desired resolution in the sense of
Stallings~\cite{Sta}, and then he introduces his
resolution as a deformation of this abstract
approximated version. Now the miraculous existence
of this top approximation directly relies on the special
symmetry of the Coxeter relations that define  Artin
monoids. For more general relations, in particular for
relations that do not preserve the length of the words,
such as those mentioned  above, even the notion of a
top factor is problematic, and extending Squier's
construction appears quite problematic---see also
Remark~\ref{R:Squi} for further comments about
obstructions to extending~\cite{Squ}.

In this paper, we develop new solutions, which address
the complete construction directly. Actually, we
propose two methods, one more simple, and one more
general. Our first solution is based on word reversing, a
syntactic technique introduced in~\cite{Dfa} for
investigating those monoids admitting least common
multiples. Starting with two words~$u$,
$v$ that represent some elements~$x$, $y$ of our
monoid,  word reversing constructs (in good cases) two
new words~$u'$, $v'$ such that both $u' v$ and $v' u$
represent the left lcm of~$x$ and~$y$, when the latter
exists. The idea here is to use word reversing to fill the
faces of the simplexes we are about to construct. The
resulting method turns out to be very simple, and we
show that it leads to a free resolution of~$\ZZ$ for
every Gaussian monoid (and even for more general
monoids called locally Gaussian) provided we start with
a convenient family of generators, typically the divisors
of the fundamental element~$\D$ in the case of a
Garside monoid. So, for instance, we obtain an explicit
resolution in the case of the braid
monoid~$B_n^+$---and of the braid
group~$B_n$---where the degree~$k$ module is
generated by the $k$-tuples of divisors of~$\D_n$. 

Our second solution is more general. It is reminiscent of
work by Kobayashi~\cite{Kob} about the homology
of rewriting systems---see also \cite{LaP, Sqw}---and
it relies on using a convenient linear ordering on the
considered generators and an induction on some
derived well-ordering of the cells. This second
construction works for arbitrary generators in all
Gaussian monoids, and, more generally, in so-called
locally left Gaussian monoids where we only assume
that any two elements that admit a common left
multiple admit a left lcm (non-spherical Artin monoids
are typical examples). The price to pay for the
generality of the construction is that we have so far no
explicit geometrical (or homotopical) interpretation for
the boundary operator and the contracting homotopy,
excepted in low degree.

With the previous tools, we reprove and extend the
results about the homology of spherical Artin groups,
and, more generally, of arbitrary Artin monoids. In
particular, we prove

\begin{thrm}\label{T:Main}
  Assume that $M$ is a finitely generated locally left
  Gaussian monoid. Then $M$ is of type~$F \! L$, in the
  sense that $\ZZ$ admits a finite free resolution
  over~$\ZZ M$.
\end{thrm}

(See Proposition~\ref{P:Main} for an explicit
bound for the length of the resolution in terms of the
cardinality of a generating set.) 

\begin{coro}
  Every Garside group~$G$ is of type~$F \! L$, \ie,
  $\ZZ$ admits a finite free resolution over~$\ZZ G$.
\end{coro}

The paper is organized as follows. In
Section~\ref{S:Gaus}, we list the needed basic
properties of (locally) Gaussian and Garside
monoids, and, in particular, we introduce word
reversing. We also recall that the homology of a monoid
satisfying Ore's embeddability conditions coincides
with the one of its group of fractions. In
Section~\ref{S:Simp}, we consider a (locally) Gaussian
monoid~$M$ and we construct an explicit 
resolution of~$\ZZ$ by a graded free
$\ZZ M$-module relying on word reversing and on the
 greedy normal form of~\cite{Eps}. We give
a natural geometrical interpretation involving 
$n$-cubes in the Cayley graph of~$M$.
In Section~\ref{S:Orde}, we consider a locally left
Gaussian monoid~$M$ (a slightly weaker hypothesis),
and we construct a second free resolution of~$\ZZ$,
relying on a well ordering of the cells. A few examples
are investigated, including the first Artin and
Birman-Ko-Lee braid monoids.

{\bf Acknowledgements.} We thank Christian Kassel for
his comments and suggestions, as well as Ruth
Charney, John Meier, and Kim Whittlesey for
interesting discussions about their independent
approach developed in~\cite{CMW}.

\section{Gaussian and Garside monoids}
\label{S:Gaus}

\subsection{Gaussian and locally Gaussian monoids}

Our notations follow those of~\cite{Squ} on the one
hand, and those of~\cite{Dfx} and~\cite{Dgk} on the
other hand. Let $M$ be a monoid. We say that $x$ is a
\emph{left divisor} (\resp a proper left divisor) of $y$
in~$M$, denoted $x \dive y$ (\resp $x
\div y$), if $y = xz$ holds for some~$z$ (\resp for
some~$z$ with $z \not= 1$). Alternatively, we say that
$y$ is a right multiple of~$x$. Right divisors and left
multiples  are defined symmetrically (but we introduce
no specific notation).

\begin{defi}
  We say that  a monoid~$M$ is {\it left Noetherian} if
  left divisibility is well founded in~$M$, \ie, there
  exists  no infinite descending sequence $x_1 \mul
  x_2 \mul \cdots$.
\end{defi}

Note that, if $M$ is a left Noetherian monoid, there
is no invertible element in~$M$ but~$1$, and,
therefore, the relation~$\div$ is a strict ordering
on~$M$ (and so is the symmetric right divisibility
relation). For $x$, $y$ in~$M$, we say that
$z$ is a {\it least common left multiple}, or left lcm,
of~$x$ and~$y$,  if $z$ is a left multiple of~$x$
and~$y$, and every common left multiple of~$x$
and~$y$ is a left multiple of~$z$. If $z$ and $z'$ are
two left lcm's for~$x$ and~$y$, then we have $z \dive
z'$ and $z'
\dive z$ by definition, hence $z = z'$ whenever $M$ is
left Noetherian. Thus, in a left Noetherian monoid, left
lcm's are unique when they exist.

\begin{defi}
  We say that a monoid~$M$ is {\it left
  Gaussian} if it is right cancellative (\ie, $zx = zy$
  implies $x = y$), left Noetherian,
  and any two elements of~$M$ admit a left lcm.
  We say that $M$ is {\it locally left Gaussian} if it
  satisfies the first two conditions above, but the third
  one is relaxed into: any two elements that admit a
  common left multiple admit a left lcm.
\end{defi}

If $M$ is a locally left Gaussian monoid, and $x$, $y$
are elements of~$M$ that admit at least one common
left multiple, we denote by~$x \lcm y$ the left lcm
of~$x$ and~$y$,  and by~$\lc(x, y)$ the unique
element~$z$  satisfying $zy= x \lcm y$; the latter is
called the {\it left complement} of~$x$ in~$y$. Thus
we have
$$\lc(x, y) \op  y = x \lcm y = \lc(y, x) \op x$$
whenever $x$ and $y$ have a common left
multiple. Observe that, if $y$ happens to be a right
divisor of~$x$, then $\lc(x, y)$ is the corresponding
quotient, \ie, we have $x =\lc(x, y) \op y$: this should
make the notation natural. It is easy to see that, in a
locally left Gaussian monoid~$M$, any two
elements~$x$, $y$ admit a right gcd, \ie, a common
right divisor~$z$ such that every common right divisor
of~$x$ and~$y$ is a right divisor of~$z$; then $M$
equipped with right division is an inf-semi-lattice with
least element~$1$.

The notion of a (locally) right Gaussian monoid is
defined symmetrically in terms of right Noetherianity,
left cancellativity and existence of right lcm's. If $M$ is
a (locally) right Gaussian monoid, and $x$, $y$ are
elements of~$M$ that admit a common right multiple,
we denote by~$\rc(x, y)$ the unique element of~$M$
such that $x \rc(x, y)$ is the right lcm of~$x$
and~$y$, and call it the {\it right complement} of~$x$
in~$y$ (we shall need no specific notation for the right
lcm in this paper).

Finally, we introduce Gaussian monoids as those
monoids satisfying the previous conditions on both
sides:

\begin{defi}
  We say that a monoid~$M$ is {\it (locally) Gaussian} 
  if it is both (locally) left Gaussian and (locally) right
  Gaussian.
\end{defi}

Roughly speaking, Gaussian monoids are those monoids
where a good theory of divisibility exist, with in
particular left and right lcm's and gcd's for every finite
family of elements. Locally Gaussian monoids are
similar, with the exception that the lcm's operations,
and, therefore, the associated complements operations,
are only partial operations. The Artin monoid
associated with an arbitrary Coxeter matrix is a typical
example of a locally Gaussian monoid~\cite{BrS}; such
an Artin monoid is Gaussian if and only if the associated
Coxeter group is finite, \ie, in the so-called spherical
case. We refer to~\cite{Pic} and~\cite{Dgp} for many
more examples of (locally) Gaussian monoids. Let us
just still mention here the Baumslag-Solitar monoid
$\Mon(a, b; ba = ab^2)$, another typical example of a
locally left Gaussian monoid that is not Gaussian, as the
elements $ab$ and~$a$ have no common left multiple.

If $M$ is a Gaussian monoid, it satisfies Ore's
conditions~\cite{ClP} and, therefore, it embeds in a
group of fractions. We say that a group~$G$ is {\it
Gaussian} if there exists at least one Gaussian
monoid~$M$ such that $G$ is the group of fractions
of~$M$. The example of Artin's braid groups~$B_n$,
which is both the group of fractions of the
monoid~$B_n^+$~\cite{Gar} and of the
Birman-Ko-Lee monoid~$\BKL_n^+$~\cite{BKL} shows
that a given Gaussian group may be the group of
fractions of several non-isomorphic Gaussian
monoids---as well as of many more monoids that need
not be Gaussian~\cite{Pic}.

\subsection{Garside and locally Garside
monoids}

In the sequel, we shall be specially interested in finitely
generated (locally) Gaussian monoids. Actually, we
shall consider a stronger condition, namely admitting a
finite generating subset that is closed under some
operations.

\begin{defi}
  We say that a monoid~$M$ is {\it (locally) Garside}
  \footnote{Garside monoids as defined above are
  called Garside monoids in~\cite{Dgk}
  and~\cite{CMW}, but they were called ``small
  Gaussian'' or ``thin Gaussian'' in previous
  papers~\cite{Dfx, Pid}, where a more restricted notion
  of a Garside monoid was also considered.}
  if it is (locally) Gaussian and it admits a finite
  generating subset~$\XX$ that is closed under left and
  right  lcm, and under left and right complements, this
  meaning that, if $x$, $y$ belong to~$\XX$ and they
  admit a common left multiple, then the left lcm~$x
  \lcm y$ and the left complement~$\lc(x, y)$, if the
  latter is not~$1$, still belong to~$\XX$, and a similar
  condition holds with right multiples.
\end{defi}

As is shown in~\cite{Dgk}, Garside
monoids may be characterized by weaker assumptions:
for instance, a sufficient condition for a Gaussian
monoid to be Garside is to admit a finite generating
subset closed under left complement. Another
equivalent condition is the existence of a Garside
element, defined as an element~$\D$ such that the
left and right divisors of~$\D$ coincide, they are finite
in number and they generate~$M$. In this case, the
family~$\DivD$ of all divisors of~$\D$ is a finite
generating set that is closed under left and right
complement, left and right lcm, and left and right gcd.
In particular, $\DivD$ equipped with the operation of
left lcm and right gcd (or of right lcm and left gcd) is a
finite lattice, with minimum~$1$ and maximum~$\D$,
and this lattice completely determines the
monoid~$M$. It is also known that every Gaussian
monoid admits a unique minimal generating family,
which implies that it admits a unique minimal Garside
element, for instance the fundamental element~$\D_n$
in the case of the monoid~$B_n^+$ of positive braids.
Let us mention that no example of a Gaussian
non-Garside monoid of finite type is known. 

Locally Garside monoids need not possess a Garside
element~$\D$ in general. Typical examples are free
monoids and, more generally, FC-type Artin
monoids~\cite{AlC}. In the cas of a free
monoid~$\XX^*$ (the set of all words over the
alphabet~$\XX$), the set~$\XX$ is a generating set
that is trivially closed under lcm and complement: any
two distinct elements~$x$,~$y$ of~$\XX$ admit no
common multiple, so $x \lcm y$ and $\lc(x, y)$
trivially belong to~$\XX$ when they exist, \ie, never.

\subsection{Identities for the complement}

In the sequel we need a convenient lcm calculus. As
already pointed out in~\cite{Dfx, Dgk}, the main object
here is not the lcm operation, but rather the derived
complement operation and the algebraic identities it
satisfies.

\begin{nota}
  For $n \ge 2$, we write $\lc(x, {y_1, \pp, y_n})$ for 
  $\lc(x, {(y_1 \lcm \cdots \lcm y_n)})$.
\end{nota}

Thus, the iterated complement operation is defined
by the equality
  \begin{equation}\label{E:lcom}
    \lc(x, {y_1, \pp, y_n})
    \op (y_1 \lcm \cdots \lcm y_n) = x \lcm y_1 \lcm
    \cdots \lcm  y_n.
  \end{equation}
Observe that \eqref{E:lcom} remains true for $n = 0$
provided we define $\lc(x, Y)$ to be~$x$ if $Y$ is the
empty sequence.

\begin{lemm}  
  The following identities hold:
  \begin{gather}
    \label{E:Iden}
    \lc(x, {y, z}) \op \lc(y, z) = \lc({(x \lcm y)}, z),\\
    \label{E:Idep}
    \lc({(\lc(x, y))}, {(\lc(z, y))}) 
     = \lc(x, {y, z})
     = \lc({(\lc(x, z))}, {(\lc(y, z))}),\\
    \label{E:Ideq}
    \lc({(xy)}, z) = \lc(x, {(\lc(z, y))}) \op \lc(y, z),\\
    \label{E:Ider}
    \lc(z, {(xy)}) = (\lc(z, y)) \lc(, x).
  \end{gather}  
\end{lemm}

\begin{proof}
  Using the associativity of the lcm, we obtain
  $$\lc(x, {y, z}) \op \lc(y, z) \op z 
  = \lc(x, {y, z}) \op (y \lcm z)
  = x \lcm (y \lcm z)
  = (x \lcm y) \lcm z
  =  \lc({(x \lcm y)}, z) \op z,$$
  and we deduce \eqref{E:Iden} by cancelling~$z$ on
  the right. The proof of \eqref{E:Idep} is similar, as
  multiplying both $\lc({(\lc(x, y))}, {(\lc(z, y))})$ and
  $\lc(x, {y, z})$ by $y \lcm z$ on the right
  gives $x \lcm y \lcm z$. Formulas~\eqref{E:Ideq}
  and~\eqref{E:Ider} are proved by expressing in various
  ways the lcm of~$xy$ and~$z$.
\end{proof}

\subsection{Word reversing}

The constructions we shall describe in
Sections~\ref{S:Simp} and, partly, \ref{S:Orde}, rely on
a word process called {\it word reversing}. It was
introduced in~\cite{Dfa}, and investigated more
systematically in Chapter~II of~\cite{Dgd}---see
also~\cite{Dgp} for further generalizations.

If $(\XX, \RRR)$ is a monoid presentation, \ie, a set of
letters plus a list of relations~$u = v$ with $u$, $v$
words over~$\XX$, we denote by $\Mon(\XX;
\RRR)$ the associated monoid, and by $\Gr(\XX;
\RRR)$ the associated group. If $u$, $v$ are words
over~$\XX$, we shall  denote by~$\cl u$ the
element of the monoid $\Mon(\XX; \RRR)$
represented by~$u$, and we write $u \equiv v$ for $\cl
u = \cl v$. We use $\XX^*$ for the
free monoid generated by~$\XX$, \ie, the set of all
words over~$\XX$; we use $\e$ for the empty word.
We also introduce $\XX\ii$ as a disjoint copy of~$\XX$
consisting of one letter~$\a\ii$ for each letter~$\a$
of~$\XX$. Finally, we say that the presentation $(\XX,
\RRR)$ is {\it positive} if all relations in~$\RRR$ have
the form $u = v$ with $u$, $v$ nonempty, and that it is
{\it complemented} if it is positive and, for each pair of
letters~$\a$, $\b$ in~$\XX$, there exists at most one
relation of the form $v \a = u \b$ in~$\RRR$, and no
relation $u \a = v \a$ with~$u \not= v$.

\begin{defi}
  Assume that $(\XX, \RRR)$ is a positive
  monoid presentation. For $w$, $w'$ words over~$\XX
  \cup \XX\ii$, we say that $w$ is {\it
$\RRR$-reversible}
  to~$w'$ (on the left) if we can transform~$w$
  to~$w'$ by iteratively deleting subwords~$u u\ii$
  where $u$ is a word over~$\XX$, and
  replacing subwords of the form~$u v\ii$
  with~${v'}\ii u'$, where $u$, $v$ are nonempty words
  over~$\XX$  and $u' v = v' u$ is one of the relations
  of~$\RRR$.
\end{defi}

For further intuition, it is important to associate with
every reversing sequence starting with a word~$w$ a
labelled planar graph defined inductively and
analogous to a van Kampen diagram: first we associate
with~$w$ a path labelled by the successive letters
of~$w$, in which the positive letters (those in~$\XX$)
are given horizontal right-oriented edges and the
negative letters (those in~$\XX\ii$) are given vertical
down-oriented edges. Then, word reversing consists in
inductively completing the diagram by using a relation
$v' u = u' v$ of~$\RRR$ (or a trivial relation $u = u$) to
close a pattern of the form
$\vcenter{\hsize 20mm
\includegraphics{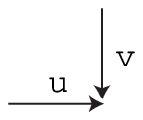}}$
 into
$\vcenter{\hsize23mm
\includegraphics{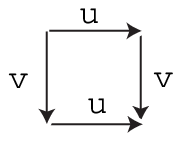}}$. 

\begin{figure}[htb]
 \includegraphics{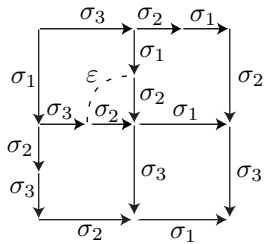}
 \caption{Word reversing diagram for 
 the standard presentation of~$B_4^+$}
  \label{F:Revg}
\end{figure}

\begin{exam}\label{X:Brai}
  Let us consider the standard presentation of the braid
  mon\-oid~$B_4^+$, namely $\Mon(\s_1, \s_2, \s_3; 
  \s_1\s_2\s_1 = \s_2\s_1\s_2, 
  \ \s_2\s_3\s_2 = \s_3\s_2\s_3, 
  \ \s_1\s_3 = \s_3\s_1)$.
  Then
  
  $\s_2\underline{\s_1\s_3\ii}\s_2\ii \to
  \underline{\s_2\s_3\ii}\s_1\s_2\ii \to
  \s_3\ii\s_2\ii\s_3\s_2\underline{\s_1\s_2\ii}
  \to$
  
  \qquad $\s_3\ii\s_2\ii\s_3\underline{\s_2\s_2\ii}
  \s_1\ii\s_2\s_1 \to
  \s_3\ii\s_2\ii\underline{\s_3\s_1\ii}\s_2\s_1 \to
  \s_3\ii\s_2\ii\s_1\ii\s_3\s_2\s_1$\\
  is a maximal reversing sequence (the pattern that is
  reversed is underlined at each step), and the
  associated diagram is displayed in Figure~\ref{F:Revg}.
\end{exam}

In general, word reversing is not a deterministic
process: starting with one word may lead to various
sequences of words, various diagrams, and, in
particular, to several terminal words, the latter being
those words that contain no pattern $\a \a\ii$ or $\a
\b\ii$ such that there exists at least one relation $v \a
= u \b$ in~$\RRR$.  However, it is easily shown that, if
$\RRR$ is a complemented presentation, then
there exists a unique maximal reversing diagram
starting with a given word~$w$, and $w$ is reversible
to at most one terminal word, so, in particular, at
most one word of the form~$u\ii v$ with
$u$, $v$ words over~$\XX$.

\begin{defi}
  Assume that $(\XX, \RRR)$ is a complemented
  presentation, and $u$, $v$ are words over~$\XX$.
  Then we denote by $\lR(u, v)$ and $\lR(v, u)$ the
  unique words over~$\XX$ such that $u v\ii$ is
  reversible to~$(\lR(v, u))\ii \lR(u, v)$, if such words
  exist.
\end{defi}

Observe that, if $\a$ and $\b$ are letters in~$\XX$,
then $\lR(\a, \b)$ and $\lR(\b, \a)$ are the (unique)
words~$u$, $v$ such that $v \a = u \b$ is a relation
in~$\RRR$, if such a relation exists.

By definition, each step of $\RRR$-reversing consists in
replacing a subword with another word that represents
the same element of the group~$\Gr(\XX; \RRR)$, so an
induction shows that, if $w$ is reversible
to~$w'$, then $w$ and
$w'$ represent the same element of~$\Gr(\XX; \RRR)$.
A slightly more careful argument gives the following
result, which is stronger in general as it need not be
true that the monoid congruence~$\equiv$ is the
restriction to positive words of the associated group
congruence, \ie, that the monoid~$\Mon(\XX; \RRR)$
embeds in the group~$\Gr(\XX; \RRR)$.

\begin{lemm}\cite{Dgp}
  Assume that $u$, $v$, $u'$, $v'$ are words
  in~$\XX^*$ and $u v\ii$ is $\RRR$-reversible
  to~${v'}\ii u'$. Then we have $\cl{v'u} =
  \cl{u' v}$, \ie, $v' u$ and $u' v$ represent
  the same element in the monoid~$\Mon(\XX; \RRR)$.
  In particular, if $(\XX, \RRR)$ is 
  complemented and $u$,~$v$ are words 
  in~$\XX^*$ such that $\lR(u, v)$ exists, we have
  $\cl{\lR(v, u)} \op \cl u = \cl{\lR(u, v)} \op \cl v$.
\end{lemm}

Thus, we see that (left) reversing constructs common
left multiples. The question is whether {\it all}
common left multiples are obtained in this way. The
answer is not always positive, but the nice point is that
there exists an effective criterion for recognizing when
this happens---and that every locally left Gaussian
monoid admits presentations for which this happens.

\begin{prop}\label{P:Pres}\cite{Dfx}
  (i) Assume that $(\XX, \RRR)$ is a complemented
  presentation satisfying the
  following conditions:

  (I) There exists a map~$\nu$ of~$\XX^*$
  to the ordinals, compatible with~$\equiv$, and
  satisfying $\nu(uv) \ge \nu(u)
  + \nu(v)$ for all~$u, v$ and $\nu(\a) > 0$ for $\a$
  in~$\XX$;   
  
  (II) We have
  $ \lR({(\lR(\a, \b))}, \, {(\lR(\g, \b))}) \equiv
  \lR({(\lR(\a, \g))}, \, {(\lR(\b, \g))})$ 
  for all $\a$, $\b$, $\g$ in~$\XX$, this meaning
  that both sides exist and are equivalent, or that neither
 exists;
  
  Then the monoid~$\Mon(\XX; \RRR)$ is
  locally left Gaussian, and, for all $u$, $v$
  in~$\XX^2$, the word~$\lR(u, v)$ exists if and only if
  the elements~$\cl u$ and $\cl v$ admit a common left
  multiple, and, in this case, $\lR(u, v)$ represents
  $\lc(\cl u, 
  \,\overline v)$; Moreover, for all words $u$, $v$, $w$,
  we have
  \begin{equation}\label{E:Cohe}
  \lR({(\lR(u, v))}, \, {(\lR(w, v))}) \equiv
  \lR({(\lR(u, w))}, \, {(\lR(v, w))}).
  \end{equation}

  (ii) Conversely, assume that $M$ is a locally left
Gaussian
  monoid, and $\XX$ is an arbitrary set of generators
  for~$M$. Let $\RRR$ consist of one relation
  $v \a = u \b$ for each pair of letters~$\a$, $\b$
  in~$\XX$ such that $\a$ and $\b$ have a common
  left multiple, where $u$ and $v$ are chosen
  (arbitrary)
   representatives of~$\lc(\a, \b)$ and $\lc(\b, \a)$
   respectively. Then $(\XX, \RRR)$ is a complemented
   presentation of~$M$ that satisfies Conditions~I
   and~II.
\end{prop}

Thus, Proposition~\ref{P:Pres} tells us that, in good
cases, left word reversing computes the left
complement operation (and, therefore, the left lcm) in
the associated monoid. If $M$ is a locally left Gaussian
monoid, and
$(\XX, \RRR)$ is a presentation of~$M$ as in
Proposition~\ref{P:Pres}(ii), then, if $\a$ and $\b$
belong to~$\XX$ and admit a common left multiple, 
the {\it word}~$\lR(\a, \b)$ of~$\XX^*$ represents the
{\it element}~$\lc(\a, \b)$ of~$M$. In particular, if
$\XX$ happens to be closed under left complement,
the word~$\lR(\a, \b)$ has length~$1$, and it consists
of the unique letter~$\lc(\a, \b)$. Thus, the
operation~$\slashstar$ can be seen as an extension of
operation~$/$ to words---as the notation suggests.
However, it should be kept in mind that
$\lR(u, v)$ is a word (not an element of the monoid),
and that computing it depends not only on~$u$, $v$,
and~$M$, but also on a particular presentation.

When $M$ is a Gaussian monoid, then, for every set of
generators~$\XX$,  Proposition~\ref{P:Pres}(ii)
provides us with a good presentation of~$M$, one
for which lcm's can be computed using word reversing.
In this case, the lcm always exists, the complement
operation is everywhere defined, and, therefore, the
operation~$\slashstar$ on words is everywhere defined
as well, which easily implies that word reversing from
an arbitrary word over~$\XX \cup \XX\ii$ always
terminates with a word~$v\ii u$ with
$u$, $v$ words over~$\XX$.

\begin{exam}
  The standard presentation of the braid
  monoid~$B_n^+$, and, more generally, the Coxeter
  presentation of all Artin monoids, are eligible for
  Proposition~\ref{P:Pres}: with a different setting, 
  verifying that Conditions~I and~II are satisfied is
  the main technical task of~\cite{Gar, BrS}, as well as
  it is the task of~\cite{BKL} in the case of the
  Birman-Ko-Lee monoid~$\BKL_n^+$.
\end{exam}

Assume that $M$ is a locally left Gaussian monoid and
$\XX$ is a generating subset of~$M$ that is closed
under left complement (a typical example is when
$M$ is a Garside monoid, and $\XX$ the set  of all
nontrivial divisors of some Garside element~$\D$ ).
Then, when applying Proposition~\ref{P:Pres}(ii), we
can choose for each pair~$\a$, $\b$ of letters, the
relation
\begin{equation}\label{E:pres}
\lc(\b, \a) \, \a = 
\lc(\a, \b) \, \b:
\end{equation}
so, here, $\lc(\a, \b)$ and $\lc(\b, \a)$ are words of
length~$1$ or $0$, \ie, letters or~$\e$. The set of these
relations, which depends only on~$M$ and on the
choice of~$\XX$, will be denoted~$\RX$ in the sequel.
As the left and the right hand sides of every relation
in~$\RX$ have length~$2$ or~$1$, $\RX$-reversing
does not increase the length of the words: for all
words~$u$, $v$ in~$\XX^*$, the length of the
word~$\lR(u, v)$ is at most the length of the
word~$u$; in particular, for every letter~$\a$ and
every word~$v$, the word~$\lR(\a, v)$ has length~$1$
or~$0$, so it is either an element of~$\XX$ or the
empty word. Another technically significant
consequence is:

\begin{lemm}\label{L:Stco}
  Assume that $M$ is a locally left Gaussian monoid,
  and $\XX$ is a generating subset of~$M$ that is 
  closed under left complement. Then the following
  strenghtening of Relation~\eqref{E:Cohe} is satisfied
  by $\RX$-reversing:
  for all words $u$, $v$, $w$ in~$\XX^*$, we have
  \begin{equation}\label{E:Stco}
  \lR({(\lR(u, v))}, \, {(\lR(w, v))}) =
  \lR({(\lR(u, w))}, \, {(\lR(v, w))}).
  \end{equation}
\end{lemm}

\begin{proof}
  Condition~II gives an equivalence for the
  words in~\eqref{E:Stco}; now, if $u$ has
  length~$1$, these words have length~$1$ at
  most, \ie, they belong to~$\XX$ or are empty, and
  equivalence implies equality for such words. The
  general case follows using an induction.
\end{proof}

\subsection{The greedy normal form}

If $M$ is a locally Gaussian monoid, and
$\XX$ is a generating subset of~$M$ that is closed
enough, we can define a unique distinguished
decomposition for every element~$x$ of~$M$ by
considering the maximal left divisor of~$x$ lying
in~$\XX$ and iterating the process. This construction
is well known in the case of Artin monoids~\cite{Dlg,
Eps, Thu, ElM}, where it is known as the (left) greedy
normal form, and it extends without change to all
Garside monoids~\cite{Dgk}. The case of locally
Gaussian monoids is not really more complicated: the
only point that could possibly fail is the
existence of a maximal divisor of~$x$ belonging
to~$\XX$; we shall see below that this existence is
guaranteed by the Noetherianity condition. Here we
describe the construction in the case of a locally {\it
right} Gaussian monoid, \ie, we use right lcm's, and not
left lcm's as in most parts of this paper:
Proposition~\ref{P:Nofo} below will explain this choice.

\begin{lemm}\label{L:Norm}
  Assume that $M$ is a locally right Gaussian monoid,
  and $\XX$ is a generating subset of~$M$ that is
  closed under right lcm. Then every nontrivial
  element~$x$ of~$M$ admits a unique greatest divisor
  lying in~$\XX$.
\end{lemm}

\begin{proof}
  Let $x = y z$ be a decomposition of~$x$ with $y \in
  \XX$ and $z$ minimal with respect to right division
  among all~$z'$ such that $x = y' z'$ holds for
  some~$y'$ in~$\XX$: such an element~$z$ exists
  since $M$ is right Noetherian. Let $y'$ be an
  arbitrary left divisor of~$x$ lying in~$\XX$. By
  construction, $y$ and $y'$ admit a common right
  multiple, namely~$x$, hence they admit a right
  lcm~$y''$ which belongs to~$\XX$, and we have $x =
  y'' z''$ for some~$z''$.  Write $y'' = y t$. Then we
  have $x = y z = y'' z'' = y t z''$, hence, by
  cancelling~$y$ on the left, $z = t z''$. The
  minimality hypothesis on~$z$ implies~$t = 1$,
  hence $y'' = y$, \ie, $y' \dive y$. So every left divisor
  of~$x$ lying in~$\XX$ is a left divisor of~$y$. 
  The uniqueness of~$y$ then follows from $1$ being
  the only invertible element of~$M$, hence the
  relation~$\dive$ being an ordering.
\end{proof}

We deduce that, under the assumptions of
Lemma~\ref{L:Norm}, every nontrivial element~$x$
of~$M$ admits a unique decomposition $x = x_1 \cdots
x_p$ such that, for each~$i$, $x_i$ is the greatest left
divisor of~
$x_i \cdots x_p$ lying in~$\XX$.  Indeed, if $x_1$
is the greatest left divisor of~$x$ lying in~$\XX$, we
have $x = x_1 x'$, and the hypothesis that $\XX$
generates~$M$ guarantees that $x_1$ is not~$1$,
hence $x'$ is a proper right divisor of~$x$, so the
hypothesis that $M$ is right Noetherian implies that
the iteration of the process terminates in a finite
number of steps.

What makes the distinguished decomposition
constructed in this way interesting is the fact that
it can be characterized using a purely local criterion,
involving only two factors at one time. This criterion is
crucial in the existence of an automatic
structure~\cite{Eps}, and it will prove crucial in our
current development as well.

\begin{defi}
  Assume that $M$ is a monoid, and $\XX$ is a subset
  of~$M$. For $x$, $y$ in~$M$, we say that $x \trX
  y$ is true if every left divisor of~$xy$ lying in~$\XX$
  is a left divisor of~$x$.
\end{defi}

\begin{lemm}\label{L:triangle}
  Assume that $M$ is a locally right Gaussian monoid,
  and $\XX$ is a generating subset of~$M$ that is
  closed under right lcm and right complement. Then
  $x \trX y \trX z$ implies $x \trX yz$.
\end{lemm}

\begin{proof}
  Let $t$ be an element of~$\XX$ dividing $xyz$ on the
  left. Let $x = x_1 \cdots x_p$ be a decomposition
  of~$x$ as a product of elements of~$\XX$. By
  hypothesis, $t$ and $x_1$ have a common right
  multiple, namely~$xyz$, hence a right lcm, say $x_1
  t_1$, and $t_1$, which is the right complement
  of~$t$ in~$x_1$, belongs to~$\XX$ by hypothesis.
  Now we have $x_1 t_1 \dive x_1 x_2 \cdots
  x_p y z$, hence $t_1 \dive x_2 \cdots
  x_p y z$. By the same argument, $t_1$ and $x_2$
  have a right lcm, say $x_2 t_2$, with $t_2 \in \XX$,
  and we have $t_2 \dive x_3 \cdots x_p y z$.
  After $p$~steps, we obtain~$t_p$ in~$\XX$ satisfying
  $t \dive x t_p$, and $t_p \dive yz$. The hypothesis $y
  \trX z$ implies $t_p \dive y$, hence $t \dive x t_p
  \dive x y$, and the hypothesis $x \trX y$ then
  implies $t \dive x$. So we proved that $t \dive xyz$
  implies $t \dive x$ for $t \in \XX$, \ie, we proved $x
  \trX xyz$.
\end{proof}

\begin{defi}
  Assume that $M$ is a monoid, and $\XX$ is a subset
  of~$M$. We say that a finite sequence $(x_1, \pp,
  x_p)$ in~$\XX^p$ is {\it $\XX$-normal} if,
  for $1 \le i < p$, we have $x_i \trX x_{i+1}$.
\end{defi}

\begin{prop}
  Assume that $M$ is a locally right Gaussian monoid,
  and $\XX$ is a generating subset of~$M$ that is
  closed under right lcm and right complement. Then
  every nontrivial element~$x$ of~$M$ admits a
  unique decomposition $x = x_1 \cdots x_p$ such
  that $(x_1, \pp, x_p)$ is a $\XX$-normal sequence.
\end{prop}

\begin{proof}
  We have already seen that every element of~$M$
  admits a unique decomposition of the form $x_1
  \cdots x_p$ with $x_1$, \pp, $x_p$ in~$\XX$
  satisfying $x_i \trX x_{i+1} \cdots x_p$ for each~$i$.
  Clearly, $x_i \trX x_{i+1} \cdots x_p$ implies
  $x_i \trX x_{i+1}$, so the only problem is to show
  that, conversely, if we have $x_1 \trX x_2 \trX \cdots
  \trX x_p$, then we have $x_i \trX x_{i+1} \cdots x_p$
  for each~$i$: this follows from
  Lemma~\ref{L:triangle} using an induction on~$p$.
\end{proof}

In the sequel, we shall denote by $\NFD x$ the
$\XX$-normal form of~$x$. For our problem, the main
property of the $\XX$-normal form is the following
connection between the normal forms of~$x$ and
of~$x \a$, established in~\cite{Dgk} in the
case of a Garside monoid:

\begin{prop}\label{P:Nofo}
  Assume that $M$ is a locally Gaussian monoid and
  $\XX$ is generating subset of~$M$ that is closed
  under right lcm, and left and right complement.
  Then, for every~$x$ in~$M$ and every~$\b$
  in~$\XX$, we have
  \begin{equation}
    \NFD x = \lR(\NFD {x \b}, \b),
  \end{equation}
  \ie, the $\XX$-normal form of~$x$ is obtained by
  reversing the word~$\NFD{x \b} \b\ii$ on the left.
\end{prop}

\begin{figure}[htb]
  $$\includegraphics{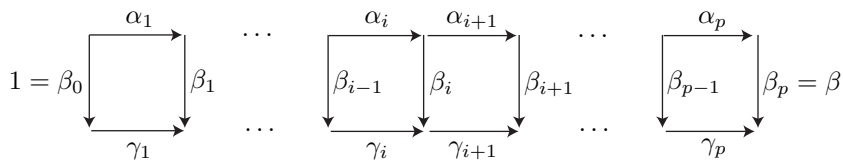}$$
  \caption{Computing the normal form using reversing}
  \label{F:Norm}
\end{figure}

\begin{proof}
  By hypothesis, the elements~$x \b$ and $\b$ admit
  a common left multiple, namely $x \b$ itself, so
  reversing the word~$\NFD{x \b} \b\ii$ on the left 
  must succeed with an empty denominator.
  Let $(\g_1, \pp, \g_p)$ be the $\XX$-normal form
  of~$x \b$. Let us define the elements~$\a_i$
  and~$\b_i$ by~$\b_p = \b$, and, using descending
  induction,
  $$\b_{i-1} = \lc(\b_i, \g_i), \qquad 
  \a_i = \lc(\g_i, \b_i)$$
  (Figure~\ref{F:Norm}).
  The hypothesis that the elements~$x \b$ and $\b$
  admit a common left multiple, namely $x \b$ itself,
  in~$M$ guarantees that $\b_i$ and $\g_i$ admit
  a common left multiple, and, therefore, the 
  inductive definition leads to no obstruction, and,
  in addition, we must have $\b_0 = 1$. By definition,
  the result of reversing $\g_1 \cdots \g_p \b\ii$ to
  the left is the word~$\a_1 \cdots \a_p$, so the
  question is to prove that $(\a_1, \pp, \a_p)$ is
  the $\XX$-normal form of~$x$. First, in~$M$, we
  have $\a_1 \cdots \a_p = \g_1 \cdots \g_p \b\ii = x
  \b \b\ii = x$, so the only question is to prove that the
  sequence $(\a_1, \pp, \a_p)$ is $\XX$-normal. 
  
  We shall prove that, for each~$i$, the relation
  $\g_i \trX \g_{i+1}$, which is true as, by
  hypothesis, the sequence $(\g_1, \pp, \g_p)$
  is $\XX$-normal, implies $\a_i \trX \a_{i+1}$.
  
  So, let us assume that some element~$\d$ of~$\XX$
  is a left divisor of~$\a_i \a_{i+1}$. Then we have
  $\d \dive \a_i \a_{i+1} \b_{i+1} = \b_{i-1} \g_i
  \g_{i-1}$. Let $\b_{i-1}\d'$ be the right lcm of~$\d$
  and~$\b_{i-1}$, which exists as $\b_{i-1} \g_i
  \g_{i+1}$ is a common right multiple of~$\d$
  and~$\b_{i-1}$. Then $\d'$ belongs to~$\XX$, and
  we have $\d' \dive \g_i \g_{i+1}$, hence $\d' \dive
  \g_i$ as $\g_i \trX g_{i+1}$ holds by hypothesis.
  Hence $\d$ is a left divisor of~$\b_{i-1} \g_i$, \ie,
  of~$\a_i \b_i$. Let $\a_i \d''$ be the right lcm
  of~$\d$ and~$\a_i$. Then $\d \dive \a_i
  \a_{i+1}$ implies $\d'' \dive \a_{i+1}$, and $\d
  \dive \a_i \b_i$ implies $\d'' \dive \b_i$. Now, by
  construction, the only common left divisor
  of~$\a_{i+1}$ and~$\b_i$ is~$1$, for, otherwise, 
  $\a_{i+1} \b_{i+1}$ would not be the left lcm
  of~$\b_{i+1}$ and~$\g_{i+1}$. So we have $\d'' =
  1$, \ie, $\d$ is a left divisor of~$\a_i$, and $\a_i
  \trX \a_{i+1}$ is true.
\end{proof}

\subsection{Group of fractions vs. monoid}

Our purpose in the sequel is to compute
the homology of a (semi)-Gaussian monoid starting
from a presentation. When the considered monoid~$M$
satisfies Ore's conditions on the left, \ie, when $M$ is
cancellative and any two elements of~$M$ admit a
common left multiple, then $M$ embeds in a group of
left fractions~$G$, and every presentation of~$M$ as a
monoid is a presentation of~$G$ as a group. By
tensorizing by~$\ZZ G$ over~$\ZZ M$ we can extend
every (left) $\ZZ M$-module into a $\ZZ G$-module.
As in in~\cite{Squ}, we shall use the following result:

\begin{prop}\cite{CaE}\label{P:Ore}
  Assume that $M$ is a monoid satisfying the Ore
  conditions on the left. Let $G$ be the group of
  fractions of~$M$. Then the functor~$R \to \ZZ G
  \otimes_{\ZZ M} R$ is exact.
\end{prop}

\begin{coro}\label{P:Grou}
  Under the above hypotheses, we have $H_*(G, \ZZ ) =
  H_*(M, \ZZ )$.
\end{coro}

So, from now on, we shall consider monoids
exclusively.  When the monoid happens to be an Ore
monoid, the homology of the monoid automatically
determines the homology of the associated group of
fractions, but the case is not really specific.

\section{The reversing resolution}
\label{S:Simp}

In this section, we assume that $M$ is a locally
Gaussian monoid, \ie, $M$ is cancellative,
left and right Noetherian, and every two elements
of~$M$ admitting a common left (\resp right) multiple
admits a left (\resp right) lcm. Next we assume that
$\XX$ is a generating subset of~$M$ that is closed
under right lcm, and under left and right complement.
Special cases are $M$ being Gaussian (in this case, lcm's
always exist), $M$ being locally Garside (in this case,
$\XX$ can be assumed to be finite), and $M$ being
Garside (both conditions simultaneously: then, we can
take for~$\XX$ the divisors of some Garside
element~$\D$).

Our aim is to construct a resolution  by free $\ZZ
M$-modules for~$\ZZ$, made into a trivial $\ZZ
M$-module by putting $x \op 1 = 1$ for every~$x$
in~$M$.

\subsection{The chain complex}

We shall consider in the sequel simplicial complexes 
associated with finite families of distinct
elements of~$\XX$ that admit a left lcm. To avoid
redundant cells, we fix a linear ordering~$<$ on~$\XX$.

\begin{defi}
  For $n \ge 0$, we denote by $\XX^{[n]}$ the family
  of all strictly increasing $n$-tuples~$(\a_1, \pp,
  \a_n)$ in~$\XX$ such that $\a_1$, \pp, $\a_n$ admit
  a left left. We denote by $\CC n$ the free $\ZZ
  M$-module generated by~$\XX^{[n]}$.  The
  generator of~$\CC n$ associated with an
  element~$A$ of~$\XX^{[n]}$ is denoted~$\c[A]$,
  and it is called an {\it $n$-cell}; the left lcm
  of~$A$ is then denoted by~$\LCM A$. The unique
  $0$-cell is denoted~$\ec$.
\end{defi}

The elements of~$\CC n$ will be called {\it
$n$-chains}. As a $\ZZ$-module, $\CC n$ is generated
by the elements of the form~$x \c[A]$ with $x \in M$;
such elements will be called {\it elementary}
$n$-chains. 

The leading idea in the sequel is to associate to each
$n$-cell an oriented $n$-cube reminiscent of a van
Kampen diagram in~$M$ and constructed using the
$\RX$-reversing process of Section~\ref{S:Gaus}. The
vertices of that cube are elements of~$M$, while the
edges are labelled by elements of~$\XX$. The $n$-cube
associated with~$\ccc[\a_1, \pp ,\a_n]$ starts from the
vertex~$1$ and ends at the vertex $\a_1 \lcm \cdots
\lcm \a_n$, so the lcm of the generators~$\a_1$, \pp,
$\a_p$ is the main diagonal of the cube, as the notation
$\LCM A$ would suggest. We start with $n$~edges
labelled~$\a_1$, \pp, $\a_n$ pointing to the final
vertex,  and we construct the other edges
backwards using left reversing, \ie, we inductively
close every pattern consisting of two converging edges
$\a$, $\b$ with two diverging edges $\lc(\b, \a)$,
$\lc(\a, \b)$. The construction terminates with~$2^n$
vertices. Finally, we associate with the elementary
$n$-chain~$x \c[A]$ the image of the $n$-cube
(associated with) $\c[A]$ under the left translation
by~$x$: the cube starts from~$x$ instead of
starting from~$1$.

\begin{figure}[htb]
 \includegraphics{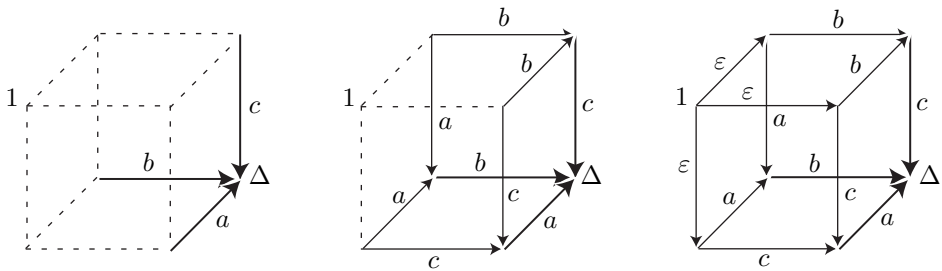}
 \caption{The $3$-cube associated with a
  $3$-cell, case of $\BKL_3^+$}\label{F:Cell}
\end{figure}

\begin{figure}[htb]
 \includegraphics{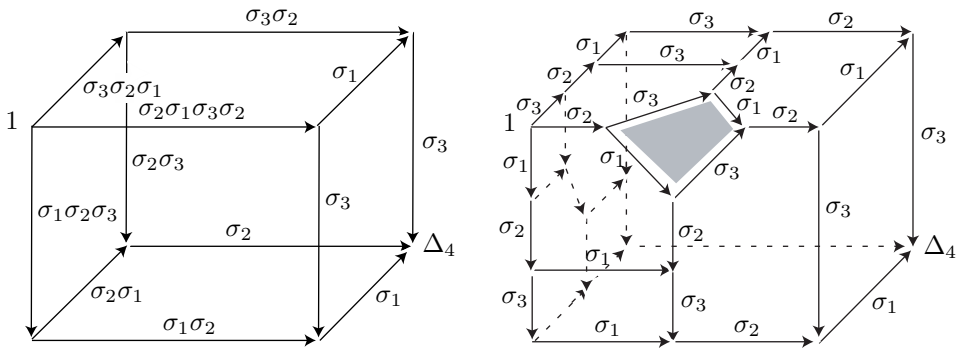}
 \caption{The $3$-cube associated with the
  $3$-cell~$\ccc[\s_1, \s_2, \s_3]$ in~$B_4^+$ when
  the generators are the divisors of~$\D_4$ (left) and
  when they are the $\s_i$'s (right)}
  \label{F:Celm}
\end{figure}

\begin{exam}\label{X:Bklt}
  Let $\BKL_3^+$ denote the Birman-Ko-Lee monoid for
  $3$-strand braids, \ie, the monoid
  $\Mon(a, b, c; ab = bc = ca)$. Then $\BKL_3^+$ is a
  Gaussian monoid, the element $\D$ defined by
  $\D= ab =bc = ca$ is a Garside element, and the
  nontrivial divisors of~$\D$ are $a$, $b$, $c$,
  and~$\D$. Thus, we can take for~$\XX$ the
  $4$-element set $\{a, b, c, \D\}$. The
  construction of the cube associated with the $3$-cell
  $\ccc[a, b, c]$ is illustrated on Figure~\ref{F:Cell}; the
  main diagonal happens to be~$\D$.
  
  Similarly, the monoid~$B_4^+$ of
  Example~\ref{X:Brai} is a Gaussian monoid, and the
  minimal Garside element is $\D_4 =
  \s_1\s_2\s_1\s_3\s_2\s_1$;  in this case, we
  can take for~$\XX$ the set of the  $23$ ($= 4! -1$)
  nontrivial divisors of~$\D_4$. The $3$-cube
  associated with the cell~$\ccc[\s_1, \s_2,
  \s_3]$ is displayed on Figure~\ref{F:Celm} (left).
\end{exam}

\begin{rema}\label{R:Atom}
  A similar construction can be made even if we
  do not assume our set of generators to be closed
  under left complement: once a complemented
  presentation has been chosen, there is a unique way
  to construct an $n$-dimensional simplex as above
  using left reversing. The only change is that words of
  length more than~$1$ appear in general, and the
  simplex, which is still a topological
  $n$-ball, is more complicated than a cube. The
  cube condition, as defined
  in~\cite{Dgk}, is the technical condition that
  guarantees that the simplex closes at the origin.
  We display on Figure~\ref{F:Celm} (right) the simplex
  associated with the $3$-cell $\ccc[\s_1, \s_2, \s_3]$
  in the standard presentation of the braid
  monoid~$B_4^+$. Observe the grey facet $\cc[\s_1,
  \s_3]$ starting at~$\s_2$: its existence corresponds to
  the fact that the {\it words} 
  $\lR({(\lR({\s_2}, {\s_1}))}, \, {(\lR({\s_3}, {\s_1}))})$
  and $\lR({(\lR({\s_2}, {\s_3}))}, \, {(\lR({\s_1}, {\s_3}))})$,
  namely $\s_2\s_1\s_3\s_2$ and
  $\s_2\s_3\s_1\s_2$,
  are equivalent, but not equal.
\end{rema}

With the previous intuition at hand, the definition of
a boundary map is clear: for $A$ an $n$-cell, we define
$\dd n \c[A]$ to be the $(n-1)$-chain obtained by
enumerating the $(n-1)$-faces of the $n$-cube
(associated with)~$\c[A]$, which are $2n$ in number,
with a sign corresponding to their orientation, and
taking into account the vertex they start from. In order
to handle such enumerations, we need to extend our
notations. 

\begin{nota}
 (i) For $\a_1$, \pp, $\a_n$ in $\XX \cup \{\e\}$, we
 define $\ccc[\a_1, \pp, \a_n]$ to be 
 $$(-1)^{\s(f)}
 \ccc[\a_{f(1)}, \pp , \a_{f(n)}]$$
  if the $\a_i$'s are pairwise distinct, $\a_{f(1)}$, \pp,
$\a_{f(n)}$ is their
  $<$-increasing enumeration, and $\s$ is the
  sign of~$f$, and to be $0_{\CC n}$
  in all other cases.

  (ii) For $A$ a cell, say $A = \ccc[\a_1, \pp
,\a_n]$,
  and $\a$ an element of~$\XX$, we denote by
  $\lc(A, \a)$ the sequence $(\lc(\a_1, \a), \pp,
  \lc(\a_n, \a))$; we denote by
  $A^i$ (\resp $A^{i, j}$) the sequence obtained by
  removing the
  $i$-th term of~$A$ (\resp  the $i$-th and the
  $j$-th terms).
\end{nota}

\begin{figure}[htb]
 \includegraphics{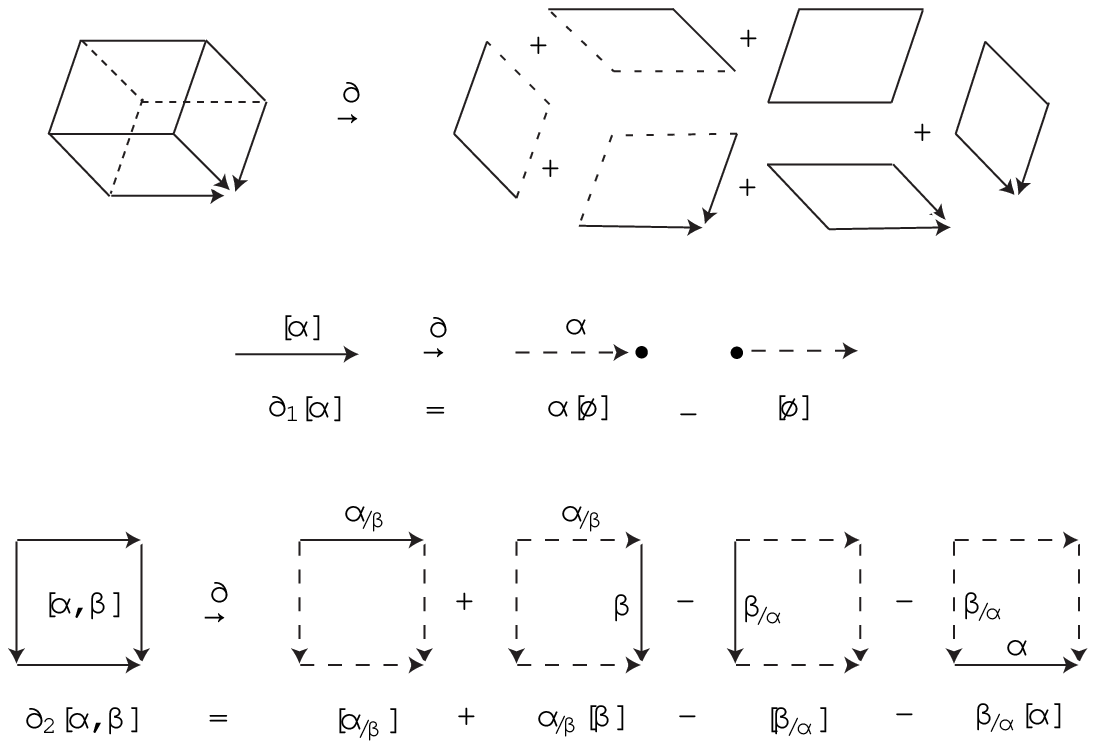}
 \caption{The boundary operator~$\dd{}$}
  \label{F:Diff}
\end{figure}

\begin{defi} (Figure~\ref{F:Diff})
  For $n \ge 1$, we define a $\ZZ M$-linear
  map $\dd n: \CC n \to \CC{n-1}$ by
  \begin{equation}\label{E:Diff}
    \dd n \c[A] = 
    \sum_{i=1}^n (-1)^i \c[\lc(A^i, \a_i)]
    - \sum_{i=1}^n (-1)^i \lc(\a_i, A^i) \c[A^i],
  \end{equation}
  for $A = (\a_1, \pp, \a_n)$; we define $\dd
  0: \CC 0 \to \ZZ$ by $\dd 0 \ec  = 1$.
\end{defi}

So, in low degrees, the formulas take the following
form:
\begin{equation}\label{E:Diff}
  \dd 1\c[\a] = \a \ec - \ec, \qquad
  \dd 2\cc[ \a , \b ] = \c[\lc(\a, \b)] + \lc(\a, \b)\c[\b] -
   \c[\lc(\b, \a)] - \lc(\b, \a)\c[\a].
\end{equation}

\begin{exam}\label{X:Bklu}
  For the Birman-Ko-Lee monoid~$\BKL3^+$,
  we read both on the above definition and on
  Figure~\ref{F:Cell} the value
  $$\dd 3 \ccc[a, b, c] = \cc[b, c] - \cc[a, c] + \cc[a,
b].$$
  Here the coefficients are $\pm 1$ as the labels of the
  three initial edges of the cube are empty words, thus
  representing~$1$ in~$M$; the three missing
  factors are $\cc[a, a]$, $\cc[b, b]$, $\cc[c, c]$, which
are null
  by definition.
  
  We suggest the reader to check on Figure~\ref{F:Celm}
  (left) the formula
  \begin{align*}
  \dd 3 \ccc[\s_1, \s_2, \s_3] = 
  - \cc[\s_1\s_2, \s_3] 
  &+ \cc[\s_2\s_1, \s_2\s_3] 
  - \cc[\s_1, \s_3\s_2] \\
  &+ \s_3\s_2\s_1 \cc[\s_2, \s_3] 
    - \s_2\s_1\s_3\s_2 \cc[\s_1, \s_3]
  + \s_1\s_2\s_3 \cc[\s_1, \s_2]
  \end{align*}
  when we consider the monoid~$B_4^+$ and
  take for~$\XX$ the divisors of the minimal Garside
  element $\D_4$.
\end{exam}

\begin{prop}
  The module $(C_*, \dd *)$ is a complex: for~$n \ge
  1$, we have $\dd
  {n-1} \dd n =~0$.
\end{prop}

\begin{proof}
  First, we have $\dd 1\c[\a]  = \a \ec  - \ec$, hence
  $\dd 0 \dd 1\c[\a] = \a \op 1 - 1_M \op 1 = 0$.
  
  Assume now $n \ge 2$. For $A = (\a_1, \pp, \a_n)$
  with $\a_1$, \pp, $\a_n \in \XX$, we obtain
  \begin{align}\label{E:deri}
    \notag \dd {n-1} \dd n \c[A] 
    & =  \sum_i (-1)^i 
    \dd{n-1}\c[\lc(A^i, \a_i)] 
     - \sum_i (-1)^i
    \lc(\a_i, A^i) \, \dd{n-1}\c[A^i]\\ 
   & =  \sum_{i \not= j}
    (-1)^{i+ j + \e(i, j)} 
    \c[\lc({(\lc(A^{i, j}, \a_i))}, {(\lc(\a_j, \a_i))})] \\
   \notag &  \hspace{1cm}
    - \sum_{i \not= j}
    (-1)^{i+ j + \e(j, i)} 
    \lc({(\lc(\a_j, \a_i))}, {(\lc(A^{i, j}, \a_i))})
   \c[\lc(A^{i, j}, \a_i)]\\
    \notag &  \hspace{1cm}
    - \sum_{i \not= j}
    (-1)^{i+ j + \e(j, i)} 
    \lc(\a_i, A^i) \, 
    \c[\lc(A^{i, j}, \a_j)] \\ 
    \notag &  \hspace{1cm}
    + \sum_{i \not= j}
    (-1)^{i+ j + \e(i, j)} 
    \lc(\a_i, A^i) \, 
    \lc(\a_j, A^{i, j}) \, 
    \c[A^{i, j}],
  \end{align}
  with $\e(i, j) = +1$ for $i < j$, and $\e(i, j) = 0$
  otherwise.  

  First, applying \eqref{E:Idep} to $\a_k$, $\a_i$, 
  and~$\a_j$, we obtain
  $\c[\lc({(\lc(A^{i, j}, \a_i))}, {(\lc(\a_j, \a_i))})]
  = \c[\lc(A^{i, j}, \a_i, \a_j)],$
  where $\a_i$ and $\a_j$ play symmetric roles, and 
  the first sum
  in~\eqref{E:deri} becomes
  $$\sum_{i \not= j} (-1)^{i+ j + \e(i, j)} 
    \c[\lc(A^{i, j}, \a_i, \a_j)].$$
  Now, each factor $\c[\lc(A^{i, j}, \a_i, \a_j)]$
  appears twice, with coefficients
  $(-1)^{i+j}$ and $(-1)^{i+j+1}$ respectively, so the
  sum vanishes.

  When applied to $\a_j$, $\a_i$, and $A^{i, j}$, 
  \eqref{E:Idep} gives
  $\lc({(\lc(\a_j, \a_i))}, {(\lc(A^{i, j}, \a_i))}) = 
  \lc(\a_j, A^j)$. It follows that the second and the
  third sum in~\eqref{E:deri}
  contain the same factors, but, as $\e(i, j) + \e(j, i) =
  1$ always holds, the signs are opposite, and the global
  sum is~$0$.
  
  Finally, applying \eqref{E:Iden} to $\a_i$, $\a_j$, and
  $\LCM{A^{i, j}\!}$ gives
   $\lc(\a_i, A^i) \,  \lc(\a_j, A^{i, j})
   = \lc({(\a_i \lcm \a_j)}, A^{i, j})$, in which $\a_i$
   and $\a_j$ play symmetric roles. So, as for
  the first sum, every factor in the fourth sum appears
  twice with opposite signs, and the sum vanishes.
  
  Observe that the case of null factors is not a problem
  above, as we always have $\lc(1, \a) = 1$ and $\lc(\a,
  1) = \a$, and, therefore, Formula~\eqref{E:Diff}
  is true for degenerate cells.
\end{proof}

It will be convenient in the sequel to extend the
notation $\ccc[\a_1, \pp, \a_n]$ to the case when the
letters~$\a_i$ are replaced by words, \ie, by finite
sequences of letters. Actually, it will be sufficient here
to consider the case when the first letter only is
replaced by a word, \ie, to consider extended cells
of the form $\cc[w, A]$ where $w$ is a word over the
alphabet~$\XX$ and $A$ is a finite sequence of
letters in~$\XX$.

\begin{defi}
  For $w$ a word over~$\XX$ and $A$ in~$\XX^{[n]}$,
  the $(n+1)$-chain $\cc[w, A]$ is defined inductively
by
  \begin{equation}\label{E:Exce}
  \cc[w, A] = 
  \begin{cases}
    0_{\CC{n+1}}
    & \mbox{if $w$ is the empty word~$\e$,} \\
    \cc[v, \lc(A, \a)] + \lc({\cl v}, {(\lc(A, \a))}) \,\,
    \cc[\a, A]     
    & \mbox{for $w = v \a$ with $\a \in \XX$.} 
  \end{cases}
  \end{equation}
\end{defi}

\begin{figure}[htb]
 \includegraphics{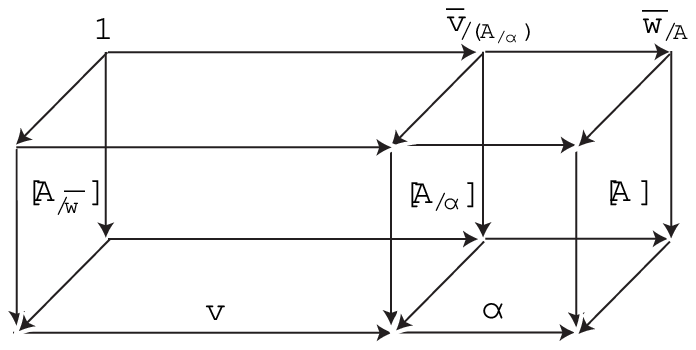}
 \caption{The chain $\cc[w, A]$ for $w = v \a$}
  \label{F:Exce}
\end{figure}

If $w$ has length~$1$, \ie, if $v$ is empty in the
inductive clause of~\eqref{E:Exce} gives $\c[v
\lc(A, \a)] = 0$ and $\lc(\cl v, {(\lc(A, \a))}) = 1$, so
our current definition of~$\cc[w, A]$ is compatible with
the previous one. Our extended notation should appear
natural when one keeps in mind the geometrical
intuition that the cell $\cc[w, A]$ is to be associated
with a $(n+1)$-parallelotope computing the left lcm
of~$\cl w$ and~$A$ using left reversing: in order to
compute the left lcm of $\cl v \a$ and $A$, we first
compute the left lcm of $\a$ and $A$, and then
compute the left lcm of~$\cl v$ and the complement
of~$A$ in~$\a$, \ie, of~$\lc(A, \a)$. However, the
rightmost cell does not start from~$1$, but from
$\lc(\cl v, {(\lc(A, \a))})$ as shown in
Figure~\ref{F:Exce}.

An easy induction shows that, for $w = \a_1\pp\a_k$,
we have for~$\c[w]$ the simple expression
\begin{equation}\label{E:Exon}
  \c[w] = \sum_{i=1}^k \cl{a_1 \cdots \a_{i-1}} \,\,
\c[\a_i].
\end{equation}
Also observe that Formula~\eqref{E:Diff} for ~$\dd 2$
can be rewritten as
\begin{equation}\label{E:Difg}
  \dd 2 \cc[ \a , \b ] = \c[\lc(\a, \b) \, \b] - \c[\lc(\b, \a) \,
\a],
\end{equation}
according to the intuition that $\dd 2$ 
enumerates the boundary of
$\vcenter{\hsize 27mm
\includegraphics{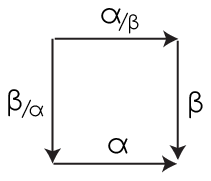}}$.

\begin{figure}[htb]
 \includegraphics{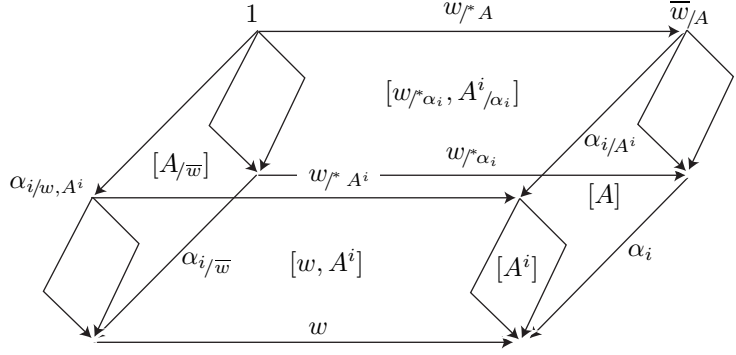}
 \caption{Decomposition of $\cc[w, A]$}
  \label{F:Excf}
\end{figure}

The following computational formula, which extends
and generalizes~\eqref{E:Diff} describes the boundary of
the parallelotope associated with~$\cc[w, A]$ taking into
account the specific role of the $w$-labelled edge: as
shown in Figure~\ref{F:Excf}, there is the right
face~$\c[A]$ at~$\lc(\cl w, A)$, the left face~$\c[\lc(A,
\,\cl w)]$, the 
$n$~lower faces $\cc[w, A^i]$ at~$\lc(\a_i,
A^i, w)$, and, finally, the $n$~upper faces
$\cc[{\lR(w, \a_i)}, {\lc(A^i, \a_i)}]$.

\begin{lemm}\label{L:exce}
  For every word~$w$, we have
  $$\dd 1 \c[w] = - \ec + \cl w \, \ec$$
  and, for $n \ge 1$ and every $n$-cell~$A$, 
  $$
      \dd{n+1} \cc[w, A] =  - \c[\lc(A, \, \cl w)]
      - \sum (-1)^i
      \cc[{\lR(w, \a_i)}, {\lc(A^i, \a_i)}]
      + \sum (-1)^i 
      \lc(\a_i, {\cl w, A^i}) \cc[w, A^i]
      + \lc(\cl w, A) \c[A] .
    $$
\end{lemm}

\begin{figure}[htb]
 \includegraphics{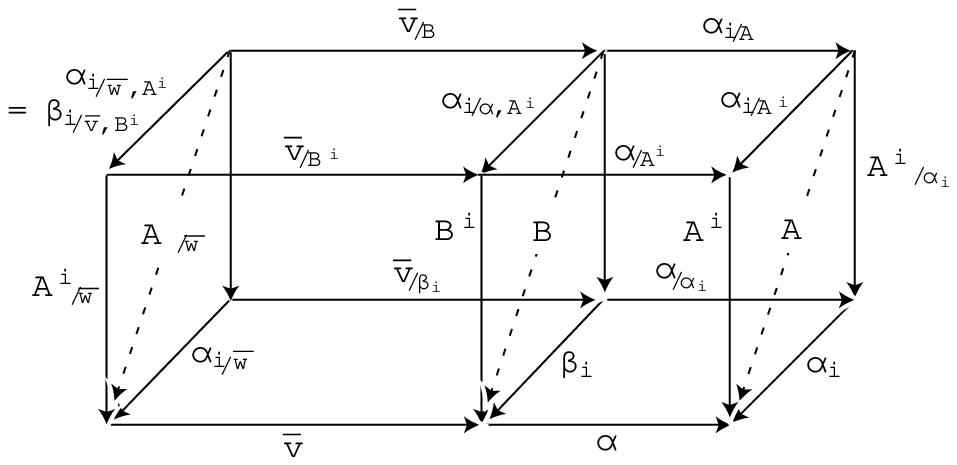}
 \caption{Computation of $\cl w \lcm \,  \LCM A$
  with $w = v \a$}
  \label{F:Expm}
\end{figure}

\begin{proof}
  The case $n = 0$ is obvious, so assume $n \ge 1$.
  We 
  use induction on the length of the word~$w$. 
  If $w$ is empty, the factors $\cc[w, A]$, $\cc[w, A^i]$,
  $\cc[{\lR(w, \a_i)}, {\lc(A^i, \a_i)}]$ vanish, 
  we have $\lc(A, \, \cl w) = A$, and the
  right hand side reduces to $\c[A] - \c[A]$, hence to~$0$,
  and the equality holds.
  Otherwise, assume $w = v \a$. By definition, we have
  $$\dd{n+1} \cc[w, A] = 
  \dd{n+1} \cc[v, B] 
   + \lc(\cl v, B) \, \dd{n+1} \cc[\a, A]
   $$
   with $B = (\b_1, \pp, \b_n) = \lc(A, \a)$.
   Applying the induction hypothesis for $\dd{n+1} \cc[v, B]$ and the definition for $\dd{n+1} \cc[\a, A]$, 
  which reads
  $$\dd{n+1}\cc[\a, A] = -[B]
  - \sum (-1)^i [\lc(\a, \a_i), \lc(A^i, \a_i)]
  + \sum (-1)^i \lc(\a_i, {\a, A^i}) [\a, A^i]
  + \lc(\a, A) [A],$$
  we obtain
  \begin{align}\label{E:comp}
        \dd{n+1} \cc[w, A] 
        & = -  \, [\lc(B, \, \cl v)]
      - \sum (-1)^i
      \cc[{\lR(v, \b_i)},  {\lc(B^i, \b_i)}]
      \\
      & \hspace{1cm} \notag
      + \sum (-1)^i
      \lc(\b_i, {\,\cl v, B^i}) 
      \cc[v, B^i]       
      + \lc(\cl v, B) \,  \c[B]
        \\
      & \hspace{1cm} \notag
      -  \lc(\cl v, B) \,  \c[B]
      - \sum (-1)^i 
      \, \lc(\cl v, B) \,  
      \cc[{\lc(\a, \a_i)}, {\lc(A^i, \a_i)}] \\
      & \hspace{1cm} \notag
      + \sum (-1)^i 
      \, \lc(\cl v, B) \,  
      \lc(\a_i, {\a, A^i}) \cc[\a, A^i]
      + \lc(\cl v, B) \,   \lc(\a, A) \c[A].
    \end{align}
   We have $\lc(\b_i, \,\cl v) =  \lc(\a_i, \,\cl w)$
  by~\eqref{E:Ider}, so the first factor
  in~\eqref{E:comp} is $-\c[\lc(A, \,\cl w)]$. Then, the
  two medial factors vanish,
   and, by construction again, we have
   $\lR(v, B) \, \lc(\a, A) = \lR(w, A)$, so
   the last factor is $\lc(\cl w, A) [A]$. There remains
   the two negative sums, and the two positive ones.
   The $i$-th factors in the negative sums are
   $$ \cc[{\lR(v, \b_i)},  {\lc(B^i, \b_i)}]
      + \lc(\cl v, B) \,  
      \cc[{\lc(\a, \a_i)}, {\lc(A^i, \a_i)}],$$
    and we claim that this is 
    $\cc[{\lR(w, \a_i)}, {\lc(A^i, \a_i)}]$. Indeed, we
have
    $\lR(w, \a_i) = \lR(v, \b_i) \, \lc(\a, \a_i)$
    as can be read on
   $\vcenter{\hsize 55mm
  \includegraphics{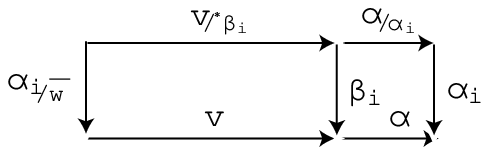}}$, so 
  \eqref{E:Exce} gives
    $$\cc[{\lR(w, \a_i)}, {\lc(A^i, \a_i)}]
    = \cc[{\lR(v, \b_i)}, 
    \lc({(\lc(A^i, \a_i))}, {(\lc(\a, \a_i))})]
    + \lc({(\lc(\cl v, \b_i))},  
    {(\lc({(\lc(A^i, \a_i))}, {(\lc(\a, \a_i))}))}) \, \, 
    \cc[{\lc(\a, \a_i)}, \lc(A^i, \a_i)].
    $$
    By~\eqref{E:Idep}, we have first
    $\lc({(\lc(A^i, \a_i))}, {(\lc(\a, \a_i))})
    = \lc({(\lc(A^i, \a))}, {(\lc(\a_i, \a))})
    = \lc(B^i, \b_i),$
    and, then, 
    $$ \lc({(\lc(\cl v, \b_i))},  
    {(\lc({(\lc(A^i, \a_i))}, {(\lc(\a, \a_i))}))})
    = \lc({(\lc(\cl v, \b_i))},  {(\lc(B^i, \b_i))})
    = \lc( \cl v, \b_i, B^i)
    = \lc( \cl v, B),$$
    which proves the claim.
    
    The argument for the positive factors
    in~\eqref{E:comp} is similar. The $i$-th factors are
    $$\lc(\b_i, {\,\cl v, B^i}) 
      \cc[v, B^i] + \lc(\, \cl v, B) \,  
      \lc(\a_i, {\a, A^i}) \cc[\a, A^i],$$
      which we claim is $\lc(\a_i, \,\cl w, A^i) \cc[w,
     A^i]$. Indeed, \eqref{E:Exce} gives
     $$\cc[w, A^i] = \cc[v, B^i] + \lc(\cl v, B^i) \, \cc[\a, 
A^i],$$
     and it remains to check the equalities
     $$\lc(\b_i, {\, \cl v, B^i})  =
     \lc(\a_i, {\, \cl w, A^i}),
     \mbox{\quad and \quad}
     \lc(\cl v, B) \op 
      \lc(\a_i, {\a, A^i}) = 
      \lc(\a_i, {\, \cl w, A^i})  \op 
      \lc(\cl v, B^i):$$
     both can be read on the diagram
    of Figure~\ref{F:Expm}, whose commutativity
  directly follows from the associativity of the lcm
  operation.
\end{proof}

\subsection{A contracting homotopy}

Our aim is to prove

\begin{prop}\label{P:main1}
  The sequence $(\CC *, \dd *)$ is a resolution
  of~$\ZZ$.
\end{prop}

To this end, it is sufficient to contruct a contracting
homotopy for~$(\CC *, \dd *)$, \ie, a family of
$\ZZ$-linear maps $\ss n: C_n \to C_{n+1}$ satisfying
$\dd{n+1} \ss n + \ss{n-1} \dd n = \id_{\CC n}$ for
each degree~$n$. We shall do it using the
$\XX$-normal form. Once again, the geometric intuition
is simple: as the chain $x \c[A]$ represents the
cube~$\c[A]$ with origin translated to~$x$, we shall
define $\ss n(x \c[A])$ to be an $(n+1)$-parallelotope
whose terminal face is~$\c[A]$ starting at~$x$. To
specify this simplex, we have to describe its $n+1$
terminal edges:  $n$ of them are the elements of~$A$;
the last one must force the main diagonal to be $x \,
\LCM A$: the most obvious choice is to take the normal
form of~$x \,\LCM A$ itself, which guarantees in
addition that there will be no initial face.

\begin{figure}[htb]
 \includegraphics{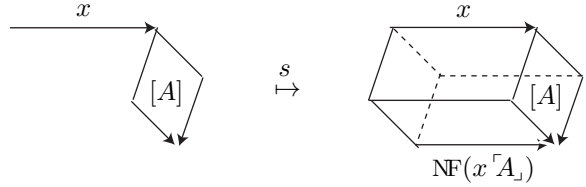}
 \caption{The contracting homotopy~$\ss{}$}
  \label{F:Sect}
\end{figure}

\begin{defi}\label{D:Sect}
  The $\ZZ$-linear mapping~$\ss n: \CC n \to
  \CC{n+1}$ is defined for $x$ in~$M$  by 
  \begin{equation}
    \ss n(x \c[A]) = \cc[\NFD{x \, \LCM A}, A]
  \end{equation}
  (Figure \ref{F:Sect});
  we define  $\ss{-1}: \ZZ \to \CC 0$ by
  $\ss{-1}(1) = \ec$.
\end{defi}

So we have in particular 
\begin{equation}\label{E:Sect}
  \ss 0(x \ec) = \c[\NFD x], \quad \mbox{and} \quad
\ss 1(x \c[\a]) = \cc[\NFD{x \a},\a]
  \end{equation}
for every~$x$ in~$M$ and every~$\a$ in~$\XX$.

\begin{lemm}
  For~$n \ge 0$, we have $\dd{n+1} \ss n +
  \ss{n-1} \dd n = \id_{\CC n}$.
\end{lemm}

\begin{proof}
  Assume first $n = 0$, and $x \in M$. Let $w = \NFD x$.
 We have
  $\ss 0(x \ec) = \c[w]$, hence $\dd 1 \ss 0(x \ec)
  = \dd 1 \c[w] = - \ec + x \ec$, and, on the other hand,
  $\dd 0(x \ec) = x \op 1 = 1$, hence
  $\ss{-1} \dd 0(x \ec) = \ec$, and 
  $(\dd 1 \ss 0 + \ss{-1} \dd 0)(x \ec) = x \ec$.
  
  Assume now $n \ge 1$. Let $w = \NFD{x \, \LCM A}$.
  Applying the definition
  of~$\ss n$ and Lemma~\ref{L:exce}, we find 
  \begin{align*}
        \dd{n+1} \ss n(x \c[A]) =  - \c[\lc(A, \,\cl w)]
      &- \sum (-1)^i
      \cc[{\lR(w, \a_i)}, \lc(A^i, \a_i)]\\
      &+ \sum (-1)^i 
      \lc(\a_i, \,\cl w, A^i) \cc[w, A^i]
      + \lc(\cl w, A) \c[A] .
  \end{align*}
  By construction, each~$\a_i$ is a right divisor of
  $\cl w$, \ie, of~$x \, \LCM A$, so we have $\c[\lc(A,
\,\cl
  w)] = \ccc[\e, \pp, \e] = 0$. At the other end, we have
  $\lc(\cl w, A) = \lc({(x \, \LCM A)}, A) = x$.
  Then $\a_i$ is a right divisor of~$\cl w$, so we have
  $\lc(\a_i, \,\cl w, A^i) = 1$, and it
  remains
  \begin{align*}
        \dd{n+1} \ss n(x \c[A]) =  
      - \sum (-1)^i
      \cc[{\lR(w, \a_i)}, \lc(A^i, \a_i)]
      + \sum (-1)^i \cc[w, A^i]
      + x \c[A] .
  \end{align*}
  
  On the other hand, we have by definition
  $$\dd n(x \c[A]) = 
    \sum_i (-1)^i x \c[\lc(A^i, \a_i)]
    - \sum_i (-1)^i x\lc(\a_i, A^i) \c[A^i].$$
  Now we have $x \LCM{\lc(A^i, \a_i)}
  \op \a_i = x \, \LCM A$, which, by
  Proposition~\ref{P:Nofo}, implies that the
  $\XX$-normal form of $x \LCM{\lc(A^i, \a_i)}$ is
  $\lR(w, \a_i)$. Then  $x\lc(\a_i, A^i)\, \LCM{A^i\!}$ is
  equal to~$x \, \LCM A$, and, therefore, its normal
form
  is~$w$. Applying the definition of~$\ss{n-1}$, we
  deduce
  $$\ss{n-1} \dd n(x \c[A]) =
  \sum (-1)^i \cc[{\lR(w, \a_i)}, \lc(A^i, \a_i)]
  - \sum (-1)^i \cc[w, A^i],$$
  and, finally, $(\dd{n+1} \ss n + \ss{n-1} \dd
  n)(x \c[A]) = x \c[A]$.
\end{proof}

Thus the sequence $\ss{*}$ is a contracting homotopy
for the complex~$(C_*, \dd{*})$, and
Proposition~\ref{P:main1} is established.

\begin{rema}
  The point in the previous argument and,
  actually, in the whole construction, is the fact that
  the normal form is computed by left reversing:
  this is what makes the explicit direct definition of
  the contracting homotopy possible. There is no need
  that the normal form we use be exactly the
  $\XX$-normal form of Section~\ref{S:Gaus}: the only
  required property is that stated in
  Proposition~\ref{P:Nofo}, namely that, if $w$ is the
  normal form of~$x \b$, then the normal form
  of~$x$ is obtained from $w$ and~$\b$ by
  left reversing.
\end{rema}

\subsection{Applications}

By definition, the set~$\XX^{[n]}$ is a basis for the
degree~$n$ module~$\CC n$ in our resolution
of~$\ZZ$ by free $\ZZ M$-modules. If the set~$\XX$
happens to be finite, then $\XX^{[n]}$ is empty for
$n$ larger than the cardinality of~$\XX$, and the
resolution is finite. By definition, choosing a
finite set~$\XX$ with the required closure properties is
possible in those monoids we called locally Garside
monoids in Section~\ref{S:Gaus}, so we may state:

\begin{prop}
  Every locally Garside monoid is of type~$F\!L$.
\end{prop}

Every Garside monoid admits a group of fractions, so,
using Proposition~\ref{P:Ore}, we deduce

\begin{coro}
  Every Garside group is of type~$F\!L$.
\end{coro}

As our constructions are explicit, they can easily be
used to compute the homology of the considered
monoid. Indeed, if we start with the resolution $(\CC *,
\dd *)$, and trivialize the elements of~$M$ in the
formulas for~$\dd *$, we obtain $\ZZ$-linear maps
that will be denoted~$\DD  *$, and we have then
$$
  H_n(M, \ZZ) = \Ker\,{\DD  n} / \Im\,{\DD  {n+1}}.
$$
Below is an example of such computations.

\begin{exam}
  Let us consider the Birman-Ko-Lee
  monoid~$\BKL_3^+$ of Example~\ref{X:Bklt} with
  $\XX = \{a, b, c, \D\}$. We recall that, by
  Proposition~\ref{X:Bklu}, the homology
  of~$\BKL_3^+$ is also that of its group of fractions,
  here the braid group~$B_3$.
  
  First, we find $\dd 1\c[a] = (a - 1) \ec$, hence $\DD 
1 \c[a]
  = 0$. The result is similar for all $1$-cells, and 
  $\Ker\,\DD  1$ is generated by $\c[a]$, $\c[b]$, $\c[c]$, 
  and~$\c[\D]$.
  
  Then, we find $\dd 2\cc[a, b] = \c[a] + a\c[b] - \c[c] - c\c[a]$,
  hence $\DD  2\cc[a, b] = \c[b] -\c[c]$, and, similarly,
  $\DD  2\cc[b, c] = \c[c] -\c[a]$, $\DD  2\cc[a, c] = \c[b] -\c[a]$, 
  $\DD  2\cc[a, \D] = \c[\D] -\c[a] - \c[c]$,
  $\DD  2\cc[b, \D] = \c[\D] -\c[b] - \c[a]$, and
  $\DD  2\cc[c, \D] = \c[\D] -\c[c] - \c[b]$. It follows that $\Im\,
  \DD  2$ is generated by the images of say $\cc[a, b]$, $\cc[a, c]$, and $\cc[b, \D]$, namely $\c[b] -\c[a]$, $\c[c] -\c[b]$,
  and $\c[\D] -\c[b] - \c[a]$, and we deduce
  $$H_1(B_3, \ZZ) = H_1(\BKL_3^+, \ZZ) =  \Ker\,\DD  1
/
  \Im\,\DD  2 = \ZZ.$$
  As we have $\DD  2\cc[b, c] = - \DD  2\cc[a, b] + \DD  2\cc[a, c]$, 
  $\DD  2\cc[a, \D] = \DD  2\cc[b, \D] + \DD  2\cc[a, b]$, and
  $\DD  2\cc[c, \D] = \DD  2\cc[a, \D] + \DD  2\cc[a, c]$,  
  $\Ker\,
  \DD  2$ is generated by
  $\cc[b, c] + \cc[a, b] - \cc[a, c]$, $\cc[a, \D] - \cc[b, \D] - \cc[a, b]$, and
  $\cc[c, \D] - \cc[a, \D] + \cc[a, c]$.
  
  Next, from the value of~$\dd 3\ccc[a, b, c]$ 
  computed in Example~\ref{X:Bklu}, we deduce
  $\DD  3\ccc[a, b, c] = \cc[b, c] - \cc[a, c] + \cc[a, b]$, and, similarly,
  $\DD  3\ccc[a, b, \D] = \cc[b, \D] - \cc[a, \D] + \cc[a, b]$,
  $\DD  3\ccc[b, c, \D] = \cc[c, \D] - \cc[b, \D] + \cc[b, c]$, and
  $\DD  3\ccc[a, c, \D] = \cc[c, \D] - \cc[a, \D] + 
  \cc[a, c]$. Therefore $\Im\,\DD  3$ is generated by
  $\cc[b, c] + \cc[a, b] - \cc[a, c]$, $\cc[a, \D] - \cc[b, \D] - \cc[a, b]$, and
  $\cc[c, \D] - \cc[a, \D] + \cc[a, c]$, 
  so it coincides with~$\Ker\, \DD  2$, and we conclude
   $$H_2(B_3, \ZZ) = H_2(\BKL_3^+, \ZZ) = \Ker\,\DD  2
  / \Im\,\DD  3 = 0.$$
  We also see that $\Ker\, \DD  3$ is generated by
  $\ccc[a, b, c] - \ccc[a, b, \D] - \ccc[b, c, \D] + \ccc[a, c, \D]$.
  
  Finally, we compute
  \begin{align*}
  \dd 4\cccc[a, b, c, \D] 
  &= - \ccc[c, c, c] + \ccc[a, a, a] - \ccc[b, b, b] +
   \ccc[\e, \e, \e] 
  + \lc(a, {b, c, \D}) \ccc[b, c, \D] \\
  & \hspace{1cm}
  - \lc(b, {a, c, \D})\ccc[a, c, \D] 
  + \lc(c, {a, b, \D})\ccc[a, b, \D] 
  - \lc(\D, {a, b, c})\ccc[a, b, c] \\
  & = \ccc[b, c, \D] - \ccc[a, c, \D] + \ccc[a, b, \D] 
  - \ccc[a, b, c]. 
  \end{align*}
  So we have $\DD  4\cccc[a, b, c, \D] = \ccc[b, c, \D] - 
  \ccc[a, c, \D]
  + \ccc[a, b, \D] - \ccc[a, b, c]$, $\Im\,\DD  4$ coincides with
  $\Ker\, \DD  3$, and $H_3(\BKL_3^+, \ZZ)$ is trivial
  (as will be obvious in the next section).
\end{exam}

\begin{rema}
  As was observed in Remark~\ref{R:Atom} and
  illustrated in Figure~\ref{F:Celm} (right), it is still
  possible to associate with every $n$-tuple of generators
  an $n$-dimensional simplex by using reversing when
  we
  consider an arbitrary set of generators~$\XX$ instead
  of the divisors of some Garside element~$\D$,
  provided Conditions~I and~II of
  Proposition~\ref{P:Pres} is satisfied. We can construct
  in this way a complex~$\CC *$, and use reversing to
  define the boundary: the formulas are not so simple
  as in \eqref{E:Diff} because the simplex is not a cube in
  general, but the principle remains the same, and a
  precise definition can be given using induction of
  $\nu(\LCM A)$, where $\nu$ is a mapping satisfying
  Condition~I.
  For instance, we obtain with the standard generators
  of~$B_4^+$ 
  \begin{align*}
    \dd 3[\s_1, \s_2, \s_3] = 
    (-1 &+ \s_1 - \s_2\s_1 + \s_3\s_2\s_1) [\s_2, \s_3]\\
     &+ (-1 +\s_2 - \s_1\s_2 - \s_3\s_2 + \s_1\s_3\s_2 -
    \s_2\s_1\s_3\s_2) [\s_1, \s_3] \\
     &+ (-1 + \s_3 - \s_2\s_3 + \s_1\s_2\s_3) [\s_1, \s_2]
  \end{align*}
  where the term $\s_2[\s_1, \s_3]$ corresponds to
  the grey facet on Figure~\ref{F:Celm} (right).
  
  The question is whether the complex is exact
  in positive degree. We have no counter-example,
  but, as no canonical normal form satisfying the
  criterion of Proposition~\ref{P:Nofo} is known in the
  general case, we do not know how to construct a
  possible contracting homotopy. As we shall
  develop below an even more general construction that
  always works, we shall leave the question open.
\end{rema}

\section{The order resolution}\label{S:Orde}

The construction of Section~\ref{S:Simp} is very simple
and convenient, but it requires using a particular set
of generators, namely one that is closed under
several operations. As a consequence, in most cases,
the simplicial complex we obtain is far from being
minimal. We shall now develop another construction,
which is more general, as it starts with an arbitrary
set of generators and does not require the
considered monoid to be locally Gaussian both on the
left and on the right. The price to pay for the extension
is that the construction of the boundary operator and of
the contracting homotopy is more complicated; in
particular, it is an inductive definition and not a direct
one as in Section~\ref{S:Simp}.

In the sequel, we assume that $M$ is a locally left
Gaussian monoid, \ie, that $M$ admits right
cancellation, that left division in~$M$ has no infinite
descending chain, and that any two elements of~$M$
that admit a common left multiple admit a left lcm. We
start with an arbitrary set of generators~$\XX$ of~$M$
that does not contain~$1$.

\subsection{Cells and chains}

Our first step is to fix a linear ordering~$<$ on~$\XX$
with the property that, for each~$x$ in~$M$, the set of
all right divisors of~$x$ is well-ordered by~$<$. At the
expense of using the axiom of choice, we can always
find such an ordering; practically, we shall be mostly
interested in the case when $\XX$ is finite, or, more
generally, when $\XX$ is possibly infinite but every
element of~$M$ can be divised by finitely many
elements of~$\XX$ only, as is the case for the direct
limit~$B_\infty^+$ of the braid monoids~$B_n^+$: in
such cases, any linear ordering on~$\XX$ is convenient.

\begin{nota}
  For $\XX$ and $<$ as above, and $x$ a
  nontrivial (\ie, not equal to~$1$) element of~$M$, we
  denote by $\md x$ the $<$-least right divisor of~$x$.
\end{nota}

As in Section~\ref{S:Simp}, the simplicial complexes
we construct are associated with finite
increasing families of generators, but we introduce
additional restrictions.

\begin{defi}
  For $n \ge 0$, we denote by~$\XXX n$ the family
  of all $n$-tuples~$(\a_1, \pp, \a_n)$ with $\a_1 <
  \cdots < \a_n \in \XX$ such that $\a_1$, \pp, $\a_n$
  admit a common left multiple (hence a left lcm), and,
  in addition, $\a_i = \md{\a_i \lcm \cdots \lcm \a_n}$
  holds for each~$i$. We let $\CCC n$ denote the free
  $\ZZ M$-module generated by~$\XXX n$.   
\end{defi}

As above, the generator of~$\CCC n$ associated with an
element~$A$ of~$\XXX n$ is denoted~$[A]$, and it
is called an {\it $n$-cell}; the left lcm of~$A$ is then
denoted by~$\LCM A$. 

\begin{exam}
  In some cases, all increasing sequences of generators
  satisfy our current additional hypotheses. For
  instance, if we consider the braid
  monoid~$B_\infty^+$ and the standard
  generators~$\s_i$ ordered by $\s_i < \s_{i+1}$, then
  there exists an $n$-cell $[\s_{i_1}, \pp, \s_{i_n}]$
  for each increasing sequence $i_1 < \pp i_n$, as
  left lcm always exist in~$B_n^+$ and $\s_{i_1}$ is
  the right divisor with least index of $\s_{i_1} \lcm
  \cdots \lcm \s_{i_n}$.
  
  On the other hand, if we consider the Birman-Ko-Lee
  monoid~$\BKL_3^+$ of Example~\ref{X:Bklt}, 
  with the ordering $a < b < c$, we see
  that there are 3 increasing sequences of length~$2$,
  namely $(a, b)$, $(a, c)$, and $(b, c)$, but there are
  two $2$-cells only, namely $\cc[a, b]$ and $\cc[a, c]$, as we have $a = \md{b \lcm c}$, which
  discards~$\cc[b, c]$.
\end{exam}

As in Section~\ref{S:Simp}, we can think of associating
with every elementary $n$-chain $x \ccc[\a_1, \pp,
\a_n]$ an $n$-dimensional oriented simplex
originating at~$x$, ending at $x(\a_1 \lcm \cdots \lcm
\a_n)$, and containing $n$ terminal edges labelled
$\a_1$, \pp, $\a_n$, but the way of filling the picture
will be different, and, in particular, the simplex is not 
a cube in general, and it seems not to be very
illuminating. The main tool here is the following
preordering on elementary chains:

\begin{defi}
  For $A$ a nonempty sequence, we denote by $\first
  A$ the first element of~$A$. Then, if $x \c[A]$, $y
\c[B]$ are elementary $n$-chains, we
  say that $x \c[A] \prec y \c[B]$ holds if we have
  either $x \, \LCM A \div y \, \LCM B$, or  $x \, \LCM A = y \, \LCM B$ and  $\first A < \first B$; for $n = 0$, we
  say that $x \ec \prec y \ec$ holds if $x \div y$ does.
  If $\sum x_i \c[A_i]$ is an arbitrary $n$-chain, we
  say that $\sum x_i \c[A_i] \prec y \c[B]$ holds if
  $x_i \c[A_i] \prec y \c[B]$ holds for every~$i$.
\end{defi}
  
\begin{lemm}\label{L:wfff}
  For every~$n$, the relation~$\prec$ on
  $n$-dimensio\-nal elementary chains is compatible
  with multiplication on the left, and it has no infinite
  decreasing sequence.
\end{lemm}

\begin{proof}
  Assume $x \c[A] \prec y \c[B]$, and let $z$ be an
  arbitrary element of~$M$. Then $x \, \LCM A \div y \, \LCM B$ implies $zx \, \LCM A \div zy \, \LCM B$, and 
  $x \, \LCM A = y \, \LCM B$ implies $zx \, \LCM A = zy \, \LCM B$,
  so we have $zx \c[A] \prec zy \c[B]$ in all cases.

  Assume now $x_1 \c[A_1] \succ x_2 \c[A_2] \succ
  \cdots$. First, we deduce $x_1 \,\LCM{A_1} \mule x_2
  \,\LCM{A_2} \mule \cdots$. As $M$ is left Noetherian,
  this decreasing sequence is eventually constant, \ie,
  for some~$i_0$, we have $x_i \, \LCM{A_i} =  x_{i+1}
  \, \LCM{A_{i+1}}$ for $i \ge i_0$. Then, for $i \ge
  i_0$,  we must have $\first{A_i} > \first{A_{i+1}}$.
  Now, by construction, $\first{A_i}$ is a right divisor
  of~$\LCM{A_i}$, hence of $x_i \, \LCM{A_i}$, and, 
  therefore, of $x_{i_0} \, \LCM{A_{i_0}}$ provided
  $i \ge i_0$ is true. But, then, the hypothesis that the
  right divisors of $x_{i_0} \, \LCM{A_{i_0}}$ are
  well-ordered by~$<$ contradicts the fact that the
  elements $\first{A_i}$ make a decreasing sequence.
\end{proof}

\subsection{Reducible chains}

We shall now construct simultaneously the boundary
maps~$\ddd n: \CCC n \to \CCC{n-1}$ together with a
contracting homotopy $\sss n :  \CCC n \to \CCC{n+1}$
and a so-called reduction map $\red n : \CCC n \to
\CCC n$. The map~$\ddd n$ is $\ZZ M$-linear,
while $\sss n$ and $\red n$ are $\ZZ$-linear.

\begin{defi}
  Assume that $x \c[A]$ is an elementary chain. We say
  that $x \c[A]$ is {\it irreducible} if either $A$ is
  empty and $x$ is~$1$, \ie, $x \, \LCM A = 1$ holds,  or
  the first element of~$A$ is the $<$-least right
  divisor of~$x \, \LCM A$, \ie, $\first A = \md{x \, \LCM A}$ holds; otherwise, we say that $x \c[A]$ is
reducible.
\end{defi}

Our construction uses induction on~$n$. The
induction hypothesis, denoted~$\PH n$, is the
conjunction of the following two statements,
where $\red n$ stands for $\sss{n-1} \comp \ddd n$: 
\begin{gather*}
\PP n \quad \ddd n( \red n (x \c[A])) =   \ddd n (x \c[A]), \\
\PQ n \quad 
\red n (x \c[A])
\begin{cases}
=  x \c[A]\mbox{ if $x \c[A]$ is irreducible,}\\
\prec  x \c[A] \mbox{ if $x \c[A]$ is reducible}
\end{cases}
\end{gather*}
(observe that $\PQ n$ makes our terminology for
reducible chains coherent). 

In degree~$0$, the construction is the same as in
Section~\ref{S:Simp}: we define $\ddd 0: \CCC 0
  \to \ZZ$ and $\sss{-1}: \ZZ \to \CCC 0$ by
  \begin{equation}
    \ddd 0(\ec) = 1, \qquad
    \sss{-1}(1) = \ec.
  \end{equation}

\begin{lemm}
  Property $\PH 0$ is satisfied.
\end{lemm}

\begin{proof}
  The mapping~$\red 0$ is $\ZZ$-linear and we have
  $$\red 0(x \ec) = \sss{-1}(\ddd 0(\ec)) = \ec$$
  for every~$x$ in~$M$. Hence, we obtain
  $$\ddd 0(\red 0(x \ec)) = \ddd 0(\ec) = 1, 
  \qquad \ddd 0(x \ec) = x \op 1 = 1$$
  owing to the trivial structure of~$\ZZ M$-module
  of~$\ZZ$. Thus $\PP 0$ holds.
  Then, by definition, $x \ec$ is irreducible if and only
  if $x$ is~$1$. In this case, we have $\red 0 (x \ec) =
  \ec$. Otherwise, we have $\red 0(x \ec) = \ec 
  \prec x \ec$ by definition of~$\red 0$, and $\PQ 0$
  holds.
\end{proof}

We assume now that $\ddd n$ and $\red n$ have been
constructed so that $\PH n$ is satisfied. We aim at
defining
$$\ddd{n+1} : \CCC{n+1} \to \CCC n, \quad
\sss n : \CCC n \to \CCC{n+1}, \quad
\red{n+1} = \sss n \comp \ddd{n+1}: \CCC{n+1} \to
\CCC{n+1}$$ 
so that $\PH{n+1}$ is satisfied.
In the sequel, we use the notation $\cc[\a, A]$
for displaying the first element of an $(n+1)$-cell; we
simply write $\LCM{a, A}$ for the associated lcm, \ie,
for $a \lcm \, \LCM A$. Thus we always have
\begin{equation}
  \LCM{\a, A} = \lc(\a, A) \op \LCM A.
\end{equation}

\begin{figure}[htb]
 \includegraphics{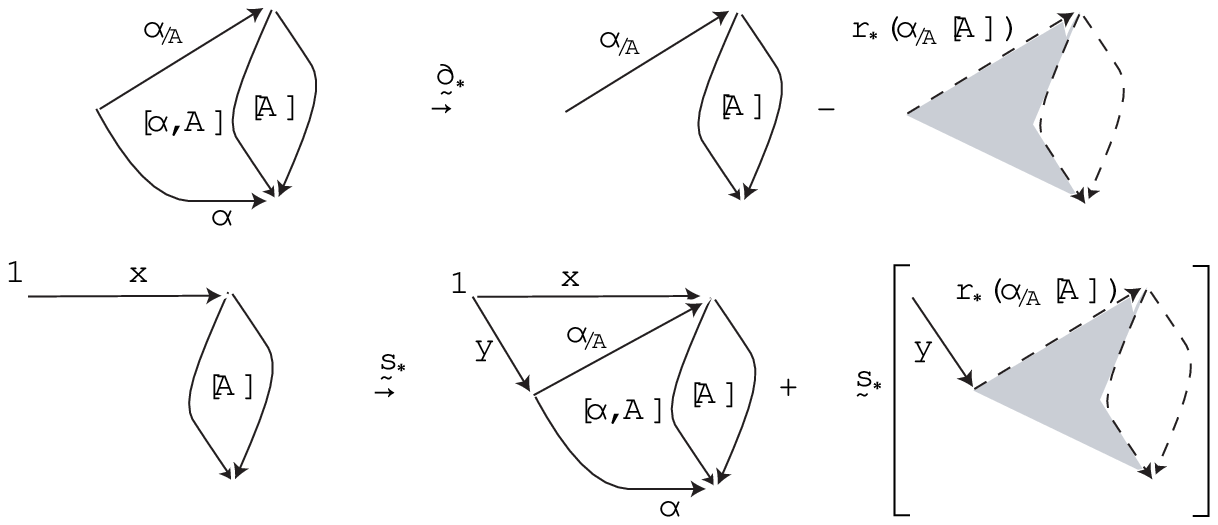}
 \caption{The boundary operator~$\ddd*$
  and the section~$\sss*$}
  \label{F:Difg}
\end{figure}

\begin{defi}
  (Figure~\ref{F:Difg})
  We define the $\ZZ M$-linear map $\ddd{n+1}:
  \CCC{n+1} \to \CCC n$ by
  \begin{equation}
    \ddd{n+1} (\cc[\a, A]) 
     =  \lc(\a, A) \c[A] - \red n(\lc(\a, A) \c[A]);
  \end{equation}
  We inductively define the $\ZZ $-linear map $\sss n:
\CCC
  n \to \CCC{n+1}$ by
  \begin{equation}\label{E:defs}
    \sss n(x \c[A]) = 
    \begin{cases}
      0 & \text{if $x \c[A]$ is irreducible},\\
      y \cc[\a, A] \rlap{+ $\sss n(y \red n (\lc(\a, A)
      [A]))$ otherwise,}\hskip1cm\\ 
      & \text{with $\a =  \md{x \, \LCM A}$ 
      and $x = y  \lc(\a, A)$.}
    \end{cases}
  \end{equation}
  Finally, we define $\red{n+1}: \CCC{n+1} \to
  \CCC{n+1}$ by $\red{n+1} = \sss n \comp
  \ddd{n+1}$.
\end{defi}

The definition of~$\ddd{n+1}$ is direct (once $\red n$
has been constructed). That of~$\sss n$ is inductive,
and we must check that it is well-founded. Now, we
observe that, in~\eqref{E:defs},  the chain~$\lc(\a, A)
\c[A]$ is reducible, as $\a < \first A$ holds by
definition, so $\PQ n$ gives $\red n(\lc(\a, A)
\c[A]) \prec \lc(\a, A) \c[A]$, and, therefore, 
  \begin{equation}\label{E:redo}
     y \red n (\lc(\a, A) \c[A]) \prec  y \lc(\a, A) \c[A]  
    = x \c[A].
  \end{equation}
Thus, our inductive definition of $\sss n$ makes sense,
and so does that of~$\red{n+1}$.

Our aim is to prove that the sequence $(\CCC *, \ddd
*)$ is a free resolution of~$\ZZ$. First, we observe that
\begin{equation}\label{E:dddd}
  \ddd n \comp \ddd{n+1} = 0
\end{equation}
automatically holds, as, using~$\PP n$, we obtain
  $$\ddd n \ddd{n+1} \cc[\a, A] 
  = \ddd n(\lc(\a, A)  \c[A]) - \ddd n (\red n
  (\lc(\a, A) \c[A])) = 0.$$

\begin{lemm}\label{L:ddse}
  Assuming $\PH n$, we have for every elementary
  $n$-chain $x \c[A]$
  \begin{equation}\label{E:Ddse}
    \ddd{n+1} \sss n (x \c[A]) 
    = x \c[A] - \red n(x \c[A]).
  \end{equation}
\end{lemm}

\begin{proof}
  We use $\prec$-induction
  on~$x \c[A]$. 
  If $x \c[A]$ is irreducible, applying $\PQ n$, we find
  $$\ddd{n+1} \sss n (x \c[A]) = 0 
  = x \c[A] - \red n (x \c[A])$$
  directly. Assume now that $x \c[A]$ is reducible. With
the
  notation of~\eqref{E:defs}, we obtain
  $$\ddd{n+1} \sss n (x \c[A])  
  = y \ddd{n+1}\cc[\a, A] 
  - \ddd{n+1} \sss n (y \red n (\lc(\a, A) \c[A])).$$
  By $\PQ n$, we have $y
  \red n (\lc(\a, A) \c[A]) \prec  x \c[A]$, so the
  induction hypothesis gives
  $$\ddd{n+1} \sss n (y \red n (\lc(\a, A) \c[A]))  
  = y \red n (\lc(\a, A) \c[A]) - \red n (y \red n
  (\lc(\a, A) \c[A])).$$
  Applying $\PP n$, we deduce
  \begin{align*}
  \red n (y \red n (\lc(\a, A) \c[A]))  
  & = \sss{n-1} (y \ddd n (\red n (\lc(\a, A) \c[A]))) \\
  & = \sss{n-1} (y \ddd n (\lc(\a, A) \c[A]))
   = \red n (y \lc(\a, A) \c[A])   
  = \red n (x \c[A]),
  \end{align*}
  hence
  $$\ddd{n+1} \sss n (x \c[A]) 
   =  y 
   \lc(\a, A) \c[A] 
   -  y \red n (\lc(\a, A) \c[A]) 
  + y \red n (\lc(\a, A) \c[A]) - \red n
  (x \c[A]),$$ 
  \ie, $\ddd{n+1} \sss n (x \c[A]) 
   =  x \c[A] - \red n (x \c[A])$, 
  as was expected.
\end{proof}

\begin{lemm}\label{L:Prpn}
  Assuming $\PH n$, $\PP{n+1}$ is satisfied.
\end{lemm}

\begin{proof}
  Assume that $x \c[A]$ is an elementary $n+1$-chain.
  We find
  \begin{align*}
    \ddd{n+1}(\red{n+1}(x \c[A]))
    &= \ddd{n+1} \sss n \ddd{n+1} (x \c[A])\\
    &= \ddd{n+1}(x \c[A]) - \red n(\ddd{n+1}(x \c[A]))\\
    &= \ddd{n+1}(x \c[A]) - \sss{n-1} \ddd n \ddd{n+1}(x
   \c[A]) = \ddd{n+1}(x \c[A]),
  \end{align*}
  by applying Lemmas~\ref{L:ddse} and~\eqref{E:dddd}.
\end{proof}

\begin{lemm}\label{L:Sedi}
  Assume that $x \cc[\a, A]$ is a reducible chain.
  Then, for each reducible chain $y \c[B]$ satisfying 
  $y \, \LCM B \dive x \, \LCM{\a, A}$, we have
  $\sss n(y \c[B]) \prec x \cc[\a, A]$.
\end{lemm}

\begin{proof}
  We use $\prec$-induction on~$y \c[B]$. By
  definition, we have
  \begin{equation}\label{E:lemm}
    \sss n(y \c[B]) = z \cc[\g, B] + \sss n(z \c[C]),
  \end{equation}
  with $\g = \md{y \, \LCM B}$, $y \, \LCM B = z 
  \, \LCM{\g, B}$, and $z \c[C]
  = z \red n(\lc(\g, B) \c[B])$.
  By \eqref{E:redo}, we always have
  $z \c[C] \prec y \c[B]$, hence, in particular, $z 
  \, \LCM C \dive y \, \LCM B \dive x \,
  \LCM{\a, A}$. So, the induction hypothesis gives
  $\sss n(z \c[C]) \prec x \cc[\a, A]$ if $z \c[C]$ is
  reducible. If $s \c[C]$ is irreducible, there is no
  contribution of $\sss n(z \c[C])$ to the sum in
  \eqref{E:lemm}, so, in both cases, it only remains to
  compare $z \cc[\g, B]$ and $x \cc[\a, A]$.

  Two cases are possible. Assume first $y \, \LCM B \div 
  x \, \LCM{\a, A}$. By construction, we have
  $z \, \LCM{\g, B} = y \, \LCM B$, so we deduce
  $z \, \LCM{\g, B} \div x \, \LCM{\a, A}$, and therefore
  $z \cc[\g, B] \prec x \cc[\a, A]$. 
  
  Assume now $y \, \LCM B = x \, \LCM{\a, A}$. By
  construction, $\g$ is the least right divisor of $y \, 
  \LCM B$, hence of $x \, \LCM{\a, A}$, and the
  hypothesis that 
  $x \cc[\a, A]$ is reducible means that $\a$ is a right
  divisor of the latter element, but is not its least right
  divisor, so we must have $\g < \a$. This,
  by definition, gives $z \cc[\g, B] \prec x \cc[\a, A]$.
\end{proof}

\begin{lemm}
  Assuming $\PH n$, $\PH{n+1}$ is satisfied.
\end{lemm}

\begin{proof}
  Owing to Lemma~\ref{L:Prpn}, it remains to prove
  $\PQ{n+1}$. Let $x \cc[\a, A]$ be an
  $(n+1)$-dimensional elementary chain. By definition,
  we have
 \begin{equation}\label{E:sedd}
   \red{n+1}(x \cc[\a, A])  
    = \sss n(x \, \lc(\a, A) \c[A]) 
    - \sss n(y \c[B]).
  \end{equation}
  with $y \c[B] = x \, \red n (\lc(\a, A) \c[A])$
  If $x \cc[\a, A]$ is irreducible, $\a$ is the least
  right divisor of $x \, \LCM{\a, A}$, the definition
  of~$\sss n$ gives
  $$\sss n(x \lc(\a, A) \c[A]) = 
  x \cc[\a, A] + \sss n(y \c[B]),$$
  and we deduce $\red{n+1}(x
  \cc[\a, A]) = x \cc[\a, A]$.
  
  Assume now that $x \cc[\a, A]$ is reducible. First, 
  we have $x \lc(\a, A) \, \LCM A = x \, \LCM{\a, A}$, so
  applying Lemma~\ref{L:Sedi} to $x \lc(\a,
  A) \c[A]$ gives $\sss n(x \lc(\a, A) \c[A]) \prec
  x \cc[\a, A]$. Then, by hypothesis, the chain $\lc(\a,
  A) \c[A]$ is reducible, so Property~$\PQ n$ gives
  $\red n(\lc(\a, A) \c[A]) \prec \lc(\a, A) \c[A]$, 
  hence, by Lemma~\ref{L:wfff}, 
  $x \, \red n(\lc(\a, A) \c[A]) \prec x \, \lc(\a, A) \c[A]$, \ie, $y \c[B] \prec x \, \lc(\a, A) \c[A]$, which
  implies in particular $z \, \LCM C \dive x \lc(\a, A)
  \, \LCM A 
  = x \,\LCM{\a, A}$. Applying Lemma~\ref{L:Sedi} to
  $y \c[B]$ gives $\sss n(y \c[B]) \prec
   x \cc[\a, A]$. Putting this in~\eqref{E:sedd}, we
  deduce $\red{n+1}(x \cc[\a, A]) \prec x \cc[\a, A]$, 
  which is Property~$\PQ{n+1}$.
\end{proof}

Thus the induction hypothesis is maintained, and the
construction can be carried out. We can now state:

\begin{prop}\label{P:main}
  The sequence $(\CCC *, \ddd *)$ is a resolution
  of~$\ZZ$ by free $\ZZ M$-modules.
\end{prop}

\begin{proof}
  First, Formula~\eqref{E:dddd} shows that $(\CCC *,
  \ddd *)$ is a complex in each degree. Then
  Formula~\eqref{E:Ddse} rewrites into
  \begin{equation}
    \ddd{n+1} \comp \sss n + \sss{n-1} \comp \ddd n =
    \id_{\CCC n},
  \end{equation}
  which shows that $\sss *$ is a contracting homotopy.
\end{proof}

An immediate corollary is the following precise version
of Theorem~\ref{T:Main}:

\begin{prop}\label{P:Main}
  Assume that $M$ is a locally left Gaussian monoid
  admitting a linearly ordered set of
  generators~$(\XX,<)$ such that $n$ is the
  maximal size of an increasing  sequence
  $(\a_1, \pp, \a_n)$ in~$\XX$ such that $\a_1 \lcm
  \cdots \lcm \a_n$ exists and $\a_i$ is the least
  right divisor of~$\a_i \lcm \cdots \lcm \a_n$ for
  each~$i$. Then $\ZZ$ admits a finite free resolution of
  length~$n$ over~$\ZZ M$; so, in particular, 
  $M$ is of type~$F\!L$.
\end{prop}

\begin{exam}
  We have seen that the Birman-Ko-Lee monoid
  $\BKL_3^+$ has a presentation with $3$ generators
  $a < b < c$, but $2$ is the maximal cardinality of a
  family as in Proposition~\ref{P:Main}, since
  $(a, b, c)$ is not eligible. We conclude 
  that $\ZZ$ admits a free resolution of
  length~$2$ over~$\ZZ\BKL_3^+$.
\end{exam}

\begin{rema}\label{R:Squi}
  Squier's approach in~\cite{Squ} has in common with
  the current approach to use the modules~$\CC n$
  (or $\CCC n$ with order assumptions dropped).
  However, the boundary operators he considers is
  different from~$\ddd n$ (and from~$\dd n$).
  Roughly speaking, Squier uses an induction on~$\div$
  and not on~$\prec$. This means that he guesses the
  exact form of all top degree factors in~$\ddd n \c[A]$, while we only guess one of these factors, namely
  the least one. Technically, the point is that, in the case
  of~\cite{Squ}, \ie, of Artin monoids, the length of the
  words induces a well defined grading on the
  monoid. Squier starts with a (very elegant)
  combinatorial construction capturing the symmetries
  of the Coxeter relations, uses it to define a first
  sketch of the differential, and then he defines his
  final differential as a deformation of the latter. It
  seems quite problematic to extend this approach to
  our general framework, because there need not
  exist any length grading, and we do not assume our
  defining relations to admit any symmetry. Due to this
  lack of symmetry, Theorem~6.10 of~\cite{Squ}, which
  is instrumental in his construction, fails in general: a
  typical example is the monoid
  $\Mon(a, b; aba= b^2)$, which is Gaussian--- the
  associated group of fractions is the braid
  group~$B_3$---and we have $\{a\} \subseteq \{a,
  b\}$, and $a \lcm b = uv$ with $u = v = b$, but
  there is no way to factor~$u = u_1u_2$, $v = v_1v_2$
  in such a way that $u_2v_1$ is equal to~$a$.
\end{rema}

\subsection{Geometrical interpretation}

We have seen that the construction of
Section~\ref{S:Simp} admits a simple geometrical
interpretation in terms of greedy normal forms and
word reversing. Here we address the question of finding
a similar geometrical interpretation for the current
construction. The answer is easy in low degree, but
quite unclear in general.

The first step is to introduce a convenient normal form
for the elements of our monoid~$M$. This is easy: as in
the case of the $\XX$-normal form, every nontrivial
element~$x$ of~$M$ has a distinguished right divisor,
namely its least right divisor~$\md x$.

\begin{defi}
 We say that a word~$w$ over~$\XX$, say $w = \a_1
 \cdots \a_p$, is the {\it ordered normal form} of~$x$,
  denoted $w = \NF x$, if  we have $x = \cl w$, and
  $\a_i = \md{\cl{\a_1 \cdots \a_i}}$ for each~$i$.
\end{defi}

Once again, an easy induction on~$\div$ shows that
every element of~$M$ admits a unique ordered normal
form: indeed, the empty word is the unique normal
form of~$1$, and, for $x \not= 1$, we write $x = y
\op \md x$, and the ordered normal form of~$x$ is
obtained by appending~$\md x$ to the
ordered normal form of~$y$.

\begin{exam}\label{X:Simp}
  Assume that $M$ is a Garside group and $\XX$ is the
  set of all divisors of some Garside element~$\D$
  of~$M$. If $<$ is any linear
  ordering on~$\XX$ that extends the opposite of
  the partial ordering given by right divisibility, then
  the ordered normal form associated with~$<$ is the
  right greedy normal form, \ie, the normal form
  constructed as the $\XX$-normal form of
  Section~\ref{S:Simp} exchanging left and right
  divisors: indeed, for every nontrivial element~$x$
  of~$M$, the rightmost factor in the right greedy
  normal form of~$x$ is the right gcd of~$x$ and~$\D$,
  hence it is a left multiple of every right divisor of~$x$
  lying in~$\XX$, and, therefore, it is the $<$-least such
  divisor.
\end{exam}

The question now is whether there exist
global expressions for $\ddd *$ and $\sss *$ in the
spirit of those of Section~\ref{S:Simp}, \ie, involving
the normal form and a word reversing process.
We still use the notation of Formula~\ref{E:Exce}, \ie, 
we write $[w]$ for the chain inductively defined
by~\eqref{E:Exce} or~\eqref{E:Exon}. 

\begin{lemm}
  For every~$x$ in~$M$ and $\a$, $\b$ in~$\XX$, we
  have
  \begin{gather}
    \ddd 1[\a] = (\a - 1) \ec, \quad
    \sss 0(x \ec) = \c[\NF x], 
    \label{E:Red0}\\
    \red 1(x \c[\a]) = \c[\NF{x \a}] - \c[\NF x], \quad
    \ddd 2 \sss 1(x \c[\a]) = \c[\NF x \, \a] -
    \c[\NF{x \a}],
    \label{E:Red1}\\
   \ddd 2 [ \a , \b ] 
   = [\NF{\lc(\a, \b)} \, \b] \!-\! [\NF{\lc(\b, \a)} \, \a]
   = [\NF{\lc(\a, \b)}] \!+\! \lc(\a, \b) [\b] 
   \!-\! [\NF{\lc(\b, \a)}] \!-\! \lc(\b, \a) [\a].
   \label{E:Difh}
  \end{gather}
\end{lemm}

\begin{proof}
  The definition gives
  $$\ddd 1[\a] = \lc(a, \emptyset) \ec - \red 0(\lc(a,
  \emptyset) \ec) = a \ec - \ec.$$
  For ~$\sss 0$, we use $\div$-induction on~$x$. For
$x = 1$, $x \ec$ is irreducible, so
   $\sss 0 (x \ec) = 0$ holds, while $\NF x$ is empty,
  and we find $\c[\NF x] = 0$. Otherwise, let $\a =
  \md x$ and $x = y a$. We have $\NF x = \NF y \op
  a$, hence $[\NF x] = [\NF y] + y [\a]$. By
  definition, we have $\sss 0(x \ec) = y [\a] + \sss
  0(y)$, hence $\sss 0(x \ec) = y [\a] + [\NF y]$ by
  induction hypothesis, and comparing the expressions
  gives $\sss 0(x \ec) = [\NF x]$.
  
  Next, we obtain
  $$\red 1(x [\a]) 
  = \sss 0(x \ddd 1[\a])
  = \sss 0(xa \ec) - \sss 0(x \ec)
  =\c[\NF{xa}] - \c[\NF{x}]$$
  The second relation in~\eqref{E:Red1} follows
  from \eqref{E:Red0} using 
  $\ddd 2 \sss 1(x [\a]) = x[\a) - \sss 0 \ddd 1(x [\a])$.
  
  Assume now that $[ \a , \b ]$ is a $2$-cell, \ie, that
  $\a < \b$ holds, $\a \lcm \b$ exists, and $\a =
  \md{\a \lcm \b}$ holds.  Applying~\eqref{E:Red1}, we
  find
  \begin{align*}\label{E:deg2}
    \ddd 2[ \a , \b ]
    &= \lc(\a, \b) [\b] - \red 1 (\lc(\a, \b) [\b]) \\
    &= \lc(\a, \b) [\b] - \c[\NF{\lc(\a, \b) \b}] 
    + [\NF{\lc(\a, \b)}]
    = [\NF{\lc(\a, \b)} \, \b] - [\NF{\a \lcm \b}].
  \end{align*}
  The hypothesis $\a = \md{\a \lcm \b}$ implies that
  the normal form of $\a \lcm \b$ is $\NF{\lc(\b, \a)}
  \a$, and we obtain
  $$\ddd 2 [ \a , \b ] 
     = [\NF{\lc(\a, \b)} \, \b] - [\NF{\lc(\b, \a)} \, \a]
     = [\NF{\lc(\a, \b)}] + \lc(\a, \b) [\b] 
     - [\NF{\lc(\b, \a)}] - \lc(\b, \a) [\a],$$
  as was expected.
\end{proof}

So, we see that the counterparts of
Formulas~\eqref{E:Diff} and~\eqref{E:Difg}, for $\ddd
1$ and $\ddd 2$ and of~\eqref{E:Sect} for~$\sss 0$
are valid: as for~$\ddd 2$, the counterpart
of~\eqref{E:Diff} has to include normal forms since, in
general, the elements~$\lc(\a, \b)$ and $\lc(\b, \a)$
do not belong to~$\XX$, as they did in the framework
of Section~\ref{S:Simp}. Observe that \eqref{E:Difh}
would fail in general if we did not restrict to cells~$[ \a , \b ]$ such that $\a$ is the least right divisor of~$\a
\lcm
\b$: this is for instance the case of the pseudo-cell~$\cc[b, c]$ in the monoid~$\BKL_3^+$ with $a < b <c$.

The next step is to interpret~$\sss 1(x[\a])$. Here, we
need to define a $2$-chain~$[u, v]$ for all word~$u$,
$v$ over~$\XX$. To this end, we keep the intuition
of Formula~\eqref{E:Exce} and use word reversing.
First, we introduce the presentation~$(\XX, \RO)$
of~$M$ by using the method of
Proposition~\ref{P:Pres}(ii) and choosing, for every
pair of letters~$\a$, $\b$ in~$\XX$, the unique
relation $\NF{\lc(\a, \b)} \, \b = \NF{\lc(\b, \a)} \,
\a$.  This presentation is uniquely determined once
$\XX$ and~$<$ have been chosen. 

\begin{defi}
  We define the $2$-chain~$\cc[u, v]$ so that the
  following rules are obeyed for all~$u$, $v$, $w$:
  $[u, \e] = 0$, $[v, u] = -[u, v]$, and
  \begin{equation}\label{E:Indu}
    [uv, w] = [u, \lo(w, v)] + 
    \lc(\cl u, {(\lc(\cl v, \cl w))})
    \: [v, w].
   \end{equation}
\end{defi}

The Noetherianity of left division in~$M$ implies that
$[u, v]$ is well defined for all~$u$,~$v$; the induction
rules mimic those of word reversing, and the idea is
that $[u, v]$ is the sum of all elementary chains
corresponding to the reversing diagram of~$u v\ii$.

\begin{ques}
  Is the following equality true:
  \begin{equation}\label{E:Seco}
  \sss 1(x \c[\a]) = \cc[\NF{x\a}, \NF x \, \a]?
  \end{equation}
\end{ques}

In the framework of Section~\ref{S:Simp}, the
definition gives
$\ss 1(x [\a]) = [\NFD{x \a}, \a]$. On the other hand, 
Formula~\ref{E:Indu} yields
$$[\NFD{x \a}, \NFD x \, \a] = 
[\NFD x, \NFD x] + \cl\e [\NFD{x \a}, \a] = 
[\NFD{x \a}, \a],$$
as, by Proposition~\ref{P:Nofo}, $\lR(\NFD{x \a}, \a)$
is equal to~$\NFD x$, and, therefore, \eqref{E:Seco}
is true. It is not hard to extend the result to our current
general framework provided the extension of
Proposition~\ref{P:Nofo} is still valid, \ie, provided
$\lo(\NF{x \a}, \a) = \NF x$ holds for every~$x$
in~$M$ and $\a$ in~$\XX$. Now, it is easy to see that
this extension is not true in general, for instance
by using the monoid~$B_4^+$ and the generators
$\s_2 < \s_1 < \s_3$. However, even if the
argument sketched above fails, Equality~\eqref{E:Seco}
remains true in all cases we tried. This suggests that
the considered geometrical interpretation could work
further.

\subsection{Examples}

Let us conclude with a few examples of our
construction. We shall successively consider the
$4$-strand braid monoid, the $3$-strand Birman-Ko-Lee
monoid, and the torus knot monoids. We use $\DDD n$
for the $\ZZ$-linear map obtained from~$\ddd n$ by
trivializing~$M$, so, again, $\Ker\, \DDD{n+1} / \Im\,
\DDD n$ is $H_n(M, \ZZ)$---as well as $H_n(G, \ZZ)$
if $M$ is a Ore monoid and $G$ is the associated group
of fractions.

\begin{exam}
Let us consider the standard presentation of~$B_4^+$.
To obtain shorter formulas, we write $a$, $b$, $c$
instead of~$\s_1$, $\s_2$, $\s_3$.  We choose $a < b
< c$. From  $\ddd 1\c[a] = (a - 1) \ec$ we deduce
$\DDD 1 [a] = 0$, and $\Ker\,\DDD 1$ is generated
by $\c[a]$, $\c[b]$, $\c[c]$. Then \eqref{E:Difh}
applies, and we find
\begin{gather*}
\ddd 2 \cc[a, b]
= \c[bab] - \c[aba]
= (-1 + b - ab) \c[a] + (1 -a + ba) \c[b], \\
\ddd 2 \cc[b, c]  
= \c[cbc] - \c[bcb]
= (-1 + c - bc) \c[b] + (1 - b + cb) \c[c], \\
\ddd 2 \cc[a, c]  
= \c[ac] - \c[ca]
= (1 - c) \c[a] + (- 1 + a) \c[c],
\end{gather*}
hence $\DDD 2( \cc[a, b]  
= - \c[a] +\c[b]$, $\DDD 2 \cc[b, c]  
= - \c[b] + \c[c]$, $\DDD 2 \cc[a, c]  = 0$.
So $\Im\, \DDD 2$ is generated by $- \c[a] +\c[b]$ and
$- \c[b] + \c[c]$, and we deduce  $H_1(B_4^+,
\ZZ) =H_1(B_4,
\ZZ) = \Ker\,\DDD 1 / \Im\,\DDD 2 = \ZZ$.

Next, we have $\LCM{a, b, c} = cbacbc = cba 
\, \LCM{b, \!c}$,
hence
$$
\ddd 3 \ccc[a, b, c]  = cba \cc[b, c] - \red 2 (cba
\cc[b, c]),
$$
and $\red 2(cba \cc[b, c])$, \ie, 
$\sss 1 \ddd 2 (cba \cc[b, c])$, evaluates to
$$
-\sss 1 (cba \c[b]) 
+ \sss 1(cbac \c[b])
- \sss 1(cbabc \c[b])
+ \sss 1 (cba \c[c]) 
- \sss 1(cbab \c[c])
+ \sss 1(cbacb \c[c]).
$$
None of the previous six chains~$x \c[A]$ is irreducible
as, in each case, $a$ is a right divisor of~$x \, \LCM A$. We
have $a = \md{cbab}$ and $cba b = c(a \lcm b)$, 
hence 
$$
\sss 1 (cba \c[b]) = c \cc[a, b] + \sss 1(c \, \red 1( ba \c[b])).
$$
Using \eqref{E:Red1}, we find
$$
\red 1 (ba \c[b])
= \c[\NF{bab}] - \c[\NF{ba}]
= \c[aba] - \c[ba]
= (1 - b + ab) \c[a] + (-1 + a) \c[b], 
$$
hence
$$
\sss 1 (cba \c[b]) = c \cc[a, b] + \sss 1(c \c[a])
- \sss 1 (cb \c[a]) + \sss 1(cab \c[a]) - \sss 1(c \c[b])
+ \sss 1(ca \c[b]).
$$
Every chain $x \c[a]$ is irreducible, and so are $c \c[b]$ and $ca \c[b]$, as we have $\md{cb} = \md{cab} = b$.
We deduce $\sss 1 (cba \c[b]) = c \cc[a, b]$. Similar
computations give $\sss 1 (cbac \c[b]) = bc \cc[a, b]$,
$\sss 1 (cbabc \c[b]) = abc \cc[a, b]$,
$\sss 1(cba \c[c]) = \cc[b, c] + cb \cc[a, c]$,
$\sss 1(cbab \c[c]) = a \cc[b, c] + (-1 + cab) \cc[a, c]$,
$\sss 1(cbacb \c[c]) = \cc[a, b] + ba \cc[b, c] + (-b + ab
+ bcab) \cc[a, c]$. We deduce the value of
$\red 2 (cba\cc[b, c])$, and, finally,
\begin{align*}
\ddd 3 \ccc[a, b, c]
&= (-1 + c - bc + abc) \cc[a, b]
+ (- 1 + a - ba + cba) \cc[b, c]\\ 
&\qquad + (-1 + b - ab - cb + cab - bcab) \cc[a, c].
\end{align*}
Triviliazing~$B_4^+$ gives
$\DDD 3 \ccc[a, b, c] = -2 \cc[a, c]$. So $\Im\,
\DDD 3$ is generated by $2 \cc[a, c]$, while $\Ker\,
\DDD 2$ is generated by $\cc[a, c]$. We deduce
$H_2(B_4^+, \ZZ) = H_2(B_4, \ZZ) = \ZZ / 2\ZZ$.

It can be observed that the values obtained above for
$\ddd *$ coincide with those of~\cite{Squ}---more
precisely, the formulas of~\cite{Squ} correspond to
what we would obtain here starting with the
initial ordering $a > b > c$: this is natural as the
presentation has the property that, for each finite
sequence of generators~$A$ in the considered
presentation, we have $\inf A = \md{\LCM A}$.
\end{exam}

\begin{exam}
  Let us consider the Birman-Ko-Lee
  monoid~$\BKL_3^+$ of Example~\ref{X:Bklt}. As the 
  minimal divisor of the lcm of~$\a$ and~$\b$ need
  not be $\a$ or $\b$, the computations are slightly
  more complicated. The reader can check that $\Ker\,
  \DDD 1$ is generated by $[a]$, $[b]$, and $[c]$,
  and that we have $\DDD 2 \cc[a, b]  = \c[b] - \c[c]$, 
  $\DDD 2\cc[b, c] = \DDD 2\cc[a, c]
  = - \c[a] + \c[b]$, so 
  $\Im\,\DDD 2$ is generated by $\c[a] - \c[b]$
  and $\c[b] - \c[c]$, and we have $H_1(M, \ZZ) = \ZZ$.
  
  For degree~2, the definition gives
  $$
  \ddd 3 \ccc[a, b, c] = \cc[b, c] - \red 2 \cc[b, c].
  $$
  Then we have
  $$
  \red 2 \cc[b, c]
  = \sss 1 \ddd 2 \cc[b, c]
  = - \sss 1 ( c \c[a]) + \sss 1 (\c[b]) - \sss 1(\c[c]) + \sss 1
  (b \c[c])
  = \sss 1 (b \c[c]),$$
  as $c \c[a]$, $\c[b]$, and $\c[c]$ are irreducible
  chains. Now we obtain
  $$
  \sss 1 (b \c[c])
  = \cc[a, c]  + \sss 1 (\red 1 (b \c[c]))
  $$
  and, by \eqref{E:Red1},
  $$
  \red 1 (b \c[c])
  = \c[\NF{bc}] - \c[\NF b]
  = \c[ca] - \c[b]
  = \c[c] + c \c[a] - \c[b],
  $$
  hence $\sss 1(\red 1 (b \c[c])) = 0$, and $\sss 1 (b 
  \c[c])
  = \cc[a, c]$. Finally, we obtain $\red 2 \cc[b, c] =
  \cc[a, c]$, and $\ddd 3 \ccc[a, b, c] = \cc[b, c] -
  \cc[a, c]$. We deduce $\DDD 3\ccc[a, b, c] = \cc[b, c]
  - \cc[a, c]$, so $\Im\, \DDD 3$ is generated by
  $\cc[b, c] - \cc[a, c]$, as is $\Ker\, \DDD 2$, and,
  therefore,  $H_2(M, \ZZ) =0$, as could be expected
  since the group of fractions of~$M$ is~$B_3$.
\end{exam}

\begin{exam}
  Finally, let $M$ be the monoid $\Mon(a, b; a^p =
b^q)$
  with $(p, q) = 1$.
  Then $M$ is locally left Gaussian, even Gaussian, and
the associated group of
  fractions is the group of the torus knot~$K_{p, q}$.
  One obtains
  $$\ddd 2 \cc[a, b] = (1 + \cdots + a^{p-1}) \c[a] +
  (1 + \cdots + b^{q-1}) \c[b],$$
  whence $\DDD 2 \, [a, b] = -p \c[a] + q \c[b]$, and 
  $H_1(M, \ZZ) = \ZZ$.
\end{exam}

\end{document}